\documentclass[12pt,reqno,twoside]{amsart}

\usepackage{amsfonts}
\usepackage{amsmath,amssymb,amsthm,amsthm}
\usepackage{mathtools}
\usepackage{dsfont}
\usepackage{etoolbox}
\usepackage{color}
\usepackage{bm}
\usepackage{subfigure}
\usepackage{units}
\usepackage{enumitem}

\usepackage{hyperref,graphicx}
\hypersetup{%
        colorlinks,
        breaklinks=true,
        plainpages=false,
        citecolor=.,
        linkcolor=.,
        urlcolor=.,
        bookmarksopen=true,%
        bookmarksnumbered=true,
        bookmarksdepth=8,
        unicode=true,%
        pdfpagelayout={SinglePage},pdffitwindow=true}%

\usepackage{cleveref}
\usepackage{booktabs}

\usepackage[dvips,bottom=1.4in,right=1in,top=1in, left=1in]{geometry}

\numberwithin{figure}{section}
\numberwithin{table}{section}

\makeatletter
\def\eqnarray{\stepcounter{equation}\let\@currentlabel=\theequation
\global\@eqnswtrue
\tabskip\@centering\let\\=\@eqncr
$$\halign to \displaywidth\bgroup\hfil\global\@eqcnt\z@
  $\displaystyle\tabskip\z@{##}$&\global\@eqcnt\@ne
  \hfil$\displaystyle{{}##{}}$\hfil
  &\global\@eqcnt\tw@ $\displaystyle{##}$\hfil
  \tabskip\@centering&\llap{##}\tabskip\z@\cr}

\def\endeqnarray{\@@eqncr\egroup
      \global\advance\c@equation\m@ne$$\global\@ignoretrue}

\usepackage{placeins}

\let\Oldsection\section
\renewcommand{\section}{\FloatBarrier\Oldsection}

 \usepackage[textsize=small]{todonotes}
\input{epsf}
\def\noi{\noindent}

\def\fvec{{\bf f}}

\def\svec{{\bf s}}

\def\uvec{{\bf u}}

\def\xvec{{\bf x}}

\def\utildevec{{\tilde {\bf u}}}

\def\Dmat{{\bf D}}

\def\Imat{{\bf I}}

\def\Kmat{{\bf K}}

\def\Mmat{{\bf M}}

\def\pmb#1{\setbox0=\hbox{$#1$}%
             \kern-.027em\copy0\kern-\wd0
             \kern+.009em\copy0\kern-\wd0
             \kern+.009em\copy0\kern-\wd0
             \kern+.009em\copy0\kern-\wd0
             \kern+.009em\copy0\kern-\wd0
             \kern+.009em\copy0\kern-\wd0
             \kern+.009em\copy0\kern-\wd0
             \kern-.045em\raise+.012em\copy0\kern-\wd0
             \kern+.009em\raise+.012em\copy0\kern-\wd0
             \kern+.009em\raise+.012em\copy0\kern-\wd0
             \kern+.009em\raise-.012em\copy0\kern-\wd0
             \kern+.009em\raise-.012em\copy0\kern-\wd0
             \kern-.018em\copy0\kern-\wd0\raise-.012em\box0}
\def\Pmb#1{\setbox0=\hbox{$#1$}%
             \kern-.033em\copy0\kern-\wd0
             \kern+.011em\copy0\kern-\wd0
             \kern+.011em\copy0\kern-\wd0
             \kern+.011em\copy0\kern-\wd0
             \kern+.011em\copy0\kern-\wd0
             \kern+.011em\copy0\kern-\wd0
             \kern+.011em\copy0\kern-\wd0
             \kern-.055em\raise+.015em\copy0\kern-\wd0
             \kern+.011em\raise+.015em\copy0\kern-\wd0
             \kern+.011em\raise+.015em\copy0\kern-\wd0
             \kern+.011em\raise-.015em\copy0\kern-\wd0
             \kern+.011em\raise-.015em\copy0\kern-\wd0
             \kern-.022em\copy0\kern-\wd0\raise-.015em\box0}
\def\alphavec{{\pmb \alpha}}

\vfill

\numberwithin{equation}{section}

\title[]{ADJOINT-BASED DETERMINATION OF WEAKNESSES IN STRUCTURES}

\author{Facundo N. Airaudo}
\address{Center for Computational Fluid Dynamics and 
Department of Physics and Astronomy, 
George Mason University, Fairfax, VA 22030, USA.}
\email{fairaudo@gmu.edu}
\author{Rainald L\"ohner}
\address{Center for Computational Fluid Dynamics and 
Department of Physics and Astronomy, 
George Mason University, Fairfax, VA 22030, USA.}
\email{rlohner@gmu.edu}
\author{Roland W\"uchner}
\address{Institut für Statik und Dynamik | Institute of Structural Analysis,
Beethovenstrasse 51, 38106 Braunschweig, Germany}
\email{r.wuechner@tu-braunschweig.de}
\author{Harbir Antil}
\address{Center for Mathematics and Artificial Intelligence (CMAI) and 
Department of Mathematical Sciences, 
George Mason University, Fairfax, VA 22030, USA.}
\email{hantil@gmu.edu}

\begin{document}

\begin{abstract}
An adjoint-based procedure to determine weaknesses, or, more generally
the material properties of structures is developed and tested. Given a
series of force and deformation/strain measurements, the material
properties are obtained by minimizing the weighted differences 
between the measured and computed values. Several examples 
with truss, plain strain and volume elements show the viability, 
accuracy and efficiency of the proposed methodology using both 
displacement and strain measurements. An important finding was 
that in order to obtain reliable, convergent results the gradient 
of the cost function has to be smoothed.
\end{abstract}

\maketitle

\section{Introduction} \label{sec:introduction}

The problem of trying to determine the material properties of
a domain from loads and measurements is common to many fields.
To mention just a few:
mining (e.g. prospecting for oil and gas), medicine (e.g. trying
to infer tissue properties), engineering (e.g. trying to determine
the existence and location of fissures, aging of structures).
\par \noi
A very pressing issue at present is the aging of concrete structures
in the developed world. Many bridges (and large buildings) were 
built with
reinforced concrete after the destruction of the second world war and the
highway system that emerged thereafter. These bridges are
now 60-70 years old, about the lifespan of concrete. Therefore, it
is imperative to know their structural integrity, which implies
determining material properties from external loads and deformations.
Damage localization is especially challenging in the case of reinforced
concrete structures due to the inhomogeneous material layout and
(mostly) very voluminous, massive structures. This motivates the development
of new damage detection techniques suitable for these applications, like
the coda wave interferometry \cite{planes2013codareview,
grabke21damagedetection, grabke22damagedetection} and their connection
to the overall
system identification to ultimately establish digital twins.
The adjoint-based technique presented here is based on displacement and
strain measurements and can thus be seen as complementary. A combination
of several sensor approaches would also appear highly promising.
Another prominent example with urgent need for damage identification
are the structures in wind generators \cite{botz2020windgenerators}. 
These massive devices
are continuously subjected to large, time-dependent forces which will
surely lead to material exhaustion and aging in 20-50 years.
\par \noi
From an abstract setting, it would seem that the task of determining 
material properties from loads and measurements is an ill-posed
problem. After all, if we think of atoms, granules or some polygonal
(e.g. finite element [FEM]) subdivision of space, the amount of data 
given resides in a space of one dimension less than the data sought. 
If we think of a cuboid domain in $d$ dimensions with $N^d$ subdivisions,
the amount of information/ data given is of $O(N^{d-1})$ while the
data sought is of $O(N^d)$.
\par \noi
Another aspect that would seem to imply that this is an ill-posed
problem is the possibility that many different spatial distributions
of material properties could yield very similar or equal deformations
under loads. That this is indeed the case for some problems is shown 
below in the examples.
\par \noi
On the other hand, the propagation of physical properties (e.g.
displacements, temperature, electrical currents, etc.) through the
domain obeys physical conservation laws, i.e. some partial differential
equations (PDEs). This implies that the material properties that can give
rise to the data measured on the boundary are restricted by these
conservation laws, i.e. are constrained. This would indicate that perhaps
the problem is not as ill-posed as initially thought. 
\par \noi
As the task of damage detection is of such importance, many techniques 
have been developed over the last decades \cite{cawley1979location, 
maia1997localization, kim2004damage, rucka2006application, 
mohan2013structural, alkayem2018structural, di2022data}.
Some of these are based on changes observed in the frequency domain 
\cite{cawley1979location, maia1997localization, mohan2013structural} 
or the time domain \cite{kim2004damage, rucka2006application}, while others
are based on changes observed in displacements or strains 
\cite{alkayem2018structural, di2022data}.
\par \noi
The procedures proposed here are also based on measured forces and
displacements/strains, but use adjoint formulations 
\cite{troltzsch2010optimal, antil2018frontiers, lohner2020determination}
and smoothing of gradients to quickly localize damaged regions.
\par
We remark that damage/weakness detection from measurements falls into 
the more general class of inverse problems where material properties 
are sought based on a desired cost functional 
\cite{borrvall2003topology, salloum2022optimization, lazarov2011filters}.

\section{Determining material properties via optimization} \label{sec:det_mat_prop_optim}

The determination of material properties (or weaknesses)
may be formulated as an 
optimization problem for the strength factor $\alpha(\xvec)$ as follows: 
Given $n$ force loadings $\fvec_i, i=1,n$ and $n$ corresponding measurements at
$m$ measuring points/locations $\xvec_j, j=1,m$ of their
respective deformations $\uvec^{md}_{ij}, \ i=1,n, \ j=1,m$ or strains
$\svec^{ms}_{ij}, \ i=1,n, \ j=1,m$, obtain the spatial distribution 
of the strength factor $\alpha$ that minimizes the cost function:

\begin{equation} \label{eq:cost}
 I(\uvec_n,\alpha) = 
  {1 \over 2} \sum_{i=1}^n \sum_{j=1}^m w^{md}_{ij} 
             ( \uvec^{md}_{ij} - \Imat^d_{ij} \cdot \uvec_i )^2
+ {1 \over 2} \sum_{i=1}^n \sum_{j=1}^m w^{ms}_{ij} 
             ( \svec^{ms}_{ij} - \Imat^s_{ij} \cdot \svec_i )^2
\end{equation}
subject to the finite element description (e.g. trusses,
beams, plates, shells, solids) of the structure \cite{zienkiewicz2005finite, simo2006computational} under 
consideration (i.e. the digital twin/system \cite{mainini2015surrogate, chinesta2020virtual}):

\begin{equation} \label{eq:forward_pde}
    \Kmat \cdot \uvec_i = \fvec_i \; , \quad i = 1, n
\end{equation}
where $w^{md}_{ij}, w^{ms}_{ij}$ are displacement and strain weights, 
$\Imat^d, \Imat^s$ interpolation matrices that are used to obtain 
the displacements and strains from the finite element mesh
at the measurement locations, and $\Kmat$ the usual stiffness matrix, 
which is obtained by assembling all the element matrices:
\begin{equation} \label{eq:strength_factor}
    \Kmat = \sum_{e=1}^{N_e} \alpha_{e} \Kmat_{e}
\end{equation}
where the strength factor $\alpha_{e}$ of the elements has already been
incorporated. We note in passing that in order to ensure that $\Kmat$
is invertible and non-degenerate $\alpha_{e} > \epsilon > 0$.

\subsection{Optimization via adjoints} \label{sec:optim_adjoints}

The objective function can be extended to the Lagrangian functional
\begin{equation} \label{eq:lagrangian_func}
    L(\uvec_n,\alpha,\utildevec_n) = I(\uvec_n,\alpha)
+ \sum_{i=1}^n \utildevec^t_i \cdot ( \Kmat \cdot \uvec_i - \fvec_i )
\end{equation}
where $\utildevec_i$ are the Lagrange multipliers (adjoints).
Variation of the Lagrangian with respect to each of the measurements
then results in:
\begin{subequations} \label{eq:optimality_cond}
\begin{eqnarray}
    &{{dL} \over {d\utildevec_i}} = \Kmat \cdot \uvec_i - \fvec_i = 0 \\
    &\sum_{j=1}^m w^{md}_{ij} (\uvec^{md}_{ij} - \Imat^d_{ij} \cdot \uvec_i)
+ \sum_{j=1}^m w^{ms}_{ij} (\svec^{ms}_{ij} - \Imat^s_{ij} \cdot \svec_i)
+ \Kmat^t \cdot \utildevec_i = 0 \\
    &{{dL} \over {d\alpha_e}}     = 
\sum_{i=1}^n \utildevec_i^t \cdot {{d\Kmat} \over {d\alpha_e}} 
                            \cdot \uvec_i =
\sum_{i=1}^n \utildevec_i^t \cdot \Kmat_e \cdot \uvec_i \, .
\end{eqnarray}
\end{subequations}

\noi
The consequences of this rearrangement are profound:
\begin{itemize}
    \item The gradient of $L$, $I$ with respect to $\alpha$ may be 
obtained by solving $n$ forward and adjoint problems; i.e.
    \item The cost for the evaluation of gradients is 
{\bf independent of the number of variables used for $\alpha$} (!).
    \item For most structural problems $\Kmat=\Kmat^t$, so if a
direct solver has been employed for the forward problem, the
cost for the evaluation of the adjoint problems is negligible;
    \item For most structural problems $\Kmat=\Kmat^t$, so if an iterative solver is employed for the forward and adjoint problems, the preconditioner can be re-utilized.
\end{itemize}

\subsection{Optimization steps} \label{sec:optim_steps}

An optimization cycle using the adjoint approach is then composed 
of the following steps:

For each force/measurement pair $i$:
\begin{enumerate}[label=\arabic*.]
    \item With current $\alpha$: solve for the deformations 
$\rightarrow \uvec_i$
    \item With current $\alpha$, $\uvec_i$ and 
$\uvec^{md}_{ij}, \svec^{md}_{ij}$: solve for the adjoints 
$\rightarrow \utildevec_i$
    \item With $\uvec_i, \utildevec_i$: obtain gradients
$\rightarrow I^i_{,\alpha}=L^i_{,\alpha}$
    \item Once all the gradients have been obtained:
    \begin{enumerate} [label=4.\arabic*.]
        \item Sum up the gradients $\rightarrow I_{,\alpha} = \sum_{i=1}^n I^i_{,\alpha}$
        \item If necessary: smooth gradients $\rightarrow I^s_{,\alpha}$
        \item Update $\alpha_{new} = \alpha_{old} - \gamma I^s_{,\alpha}$.
    \end{enumerate} 
\end{enumerate}
Here $\gamma$ is a small stepsize that can be adjusted so as to obtain
optimal convergence (e.g. via a steepest descent method).

\section{Interpolation of displacements and strains} \label{sec:interpolation}

The location of a displacement or strain gauge may not coincide
with any of the nodes of the finite element mesh. Therefore, 
in general, the displacement $\uvec_i$ at a measurement location 
$\xvec^m_i$ needs to be obtained via the interpolation matrix 
$\Imat^d_i$ as follows:
\begin{equation} \label{eq:interpolation}
    \uvec_i(\xvec^m_i) = \Imat^d_i(\xvec^m_i) \cdot \uvec
\end{equation}
where $\uvec$ are the values of the deformation vector at all 
grid points.

In many cases it is much simpler to install strain gauges instead of
displacement gauges. In this case, the strains need to be obtained
from the displacement field. This can be written formally as:
\begin{equation} \label{eq:strain_deformation}
    \svec = \Dmat \cdot \uvec
\end{equation}
where the `derivative matrix' $\Dmat$ contains the local values of the
derivatives of the shape-functions of $\uvec$.
The strain at an arbitrary position $\xvec^m_i$ is obtained via 
the interpolation matrix $\Imat^s_i$ as follows:
\begin{equation} \label{eq:strain}
    \svec_i(\xvec^m_i) = \Imat^s_i(\xvec^m_i) \cdot \svec 
                      = \Imat^s_i(\xvec^m_i) \cdot \Dmat \cdot \uvec \, .
\end{equation}
Note that in many cases the strains will only be defined in the 
elements, so that the interpolation matrices for displacements 
and strains may differ.

\section{Choice of weights} \label{sec:choice_of_weights}

The cost function is given by equation \eqref{eq:cost} repeated here for clarity:
\begin{equation}
    I(\uvec_n,\alpha) = 
  {1 \over 2} \sum_{i=1}^n \sum_{j=1}^m w^{md}_{ij} 
             ( \uvec^{md}_{ij} - \Imat^d_{ij} \cdot \uvec_i )^2
+ {1 \over 2} \sum_{i=1}^n \sum_{j=1}^m w^{ms}_{ij} 
             ( \svec^{ms}_{ij} - \Imat^s_{ij} \cdot \svec_i )^2 \, .
\end{equation}
One can immediately see that the dimensions of displacements and
strains are different. This implies that the weights should be chosen
in order that all the dimensions coincide. The simplest way of
achieving this is by making the cost function dimensionless. This
implies that the displacement weights $w^{md}_{ij}$ should be of 
dimension [1/(displacement*displacement)] and the strains weights
$w^{ms}_{ij}$ should be of dimension [1/(strain*strain)]. 
Several options are possible:

\subsubsection{Local Weighting}

In this case
\begin{equation} \label{eq:local_weighting}
    w^{md}_{ij}={1 \over{(\uvec^{md}_{ij})^2}} ~~;~~ 
   w^{ms}_{ij}={1 \over{(\svec^{ms}_{ij})^2}} ~~;
\end{equation}
this works well, but may lead to an `over-emphasis' of small 
displacements/strains that are in regions of marginal interest.

\subsubsection{Average Weighting}

In this case one first obtains the average of the absolute value of the displacements/strains for a loadcase and uses
them for the weights, i.e.:
\begin{equation} \label{eq:average_weighting}
    u_{av}={{\sum_{j=1}^m |\uvec^{md}_{ij}|} \over m} ~~; ~~
   w^{md}_{ij}={1 \over u^2_{av}} ~~;
   s_{av}={{\sum_{j=1}^m |\svec^{ms}_{ij}|} \over m} ~~; ~~
   w^{ms}_{ij}={1 \over s^2_{av}} ~~;
\end{equation}
this works well, but may lead to an `under-emphasis' of small 
displacements/strains that may occur in important regions;

\subsubsection{Max weighting}

In this case one first obtains the maximum of the absolute value of the displacements/strains for a loadcase and uses
them for the weights, i.e.:
\begin{align}
\begin{split} \label{eq:max_weighting}
    u_{max}&=max(|\uvec^{md}_{ij}|, j=1,m) ~~;~~
   w^{md}_{ij}={1 \over u^2_{max}} ~~; \\
   s_{max}&=max(|\svec^{ms}_{ij}|, j=1,m) ~~; ~~
  w^{ms}_{ij}={1 \over s^2_{max}} ~~;
\end{split}
\end{align}
\noi
this also works well for many cases, but may lead to an 
`under-emphasis' of smaller displacements/strains that can occur in 
important regions;

\subsubsection{Local/Max Weighting}

In this case
\begin{equation} \label{eq:local_max_weighting}
    w^{md}_{ij}={1 \over {max(\epsilon u_{max}, |\uvec^{md}_{ij}|))^2}} ~~;~~ 
   w^{ms}_{ij}={1 \over {max(\epsilon s_{max}, |\svec^{ms}_{ij}|))^2}} ~~;
\end{equation}

$$ w^{md}_{ij}={1 \over {max(\epsilon u_{max}, |\uvec^{md}_{ij}|))^2}} ~~;~~ 
   w^{ms}_{ij}={1 \over {max(\epsilon s_{max}, |\svec^{ms}_{ij}|))^2}} ~~;
                                                       \eqno(4.5) $$

\noi
with $\epsilon=O(0.01-0.10)$; this seemed to work best of all, 
as it combines local weighting with a max-bound minimum for local values.

\section{Smoothing of gradients} \label{sec:smoothing}

The gradients of the cost function with respect to $\alpha$ allow for
oscillatory solutions. One must therefore smooth or `regularize' the
spatial distribution. This happens naturally when using few degrees of
freedom, i.e. when $\alpha$ is defined via other spatial shape functions
(e.g. larger spatial regions of piecewise constant $\alpha$). As the
(possibly oscillatory) gradients obtained in the (many) finite elements
are averaged over spatial regions, an intrinsic smoothing occurs.
This is not the case if $\alpha$ and the gradient are defined and 
evaluated in each element separately, allowing for the largest degrees of
freedom in a mesh and hence the most accurate representation.
Three different types of smoothing or `regularization' were considered.
All of them start by performing a volume averaging from elements
to points:
\begin{equation} \label{eq:averaging}
    \alpha_p = {{ \sum_{e} \alpha_{e} V_{e} } \over
               { \sum_{e} V_{e} }}
\end{equation}
where $\alpha_p, \alpha_{e}, V_{e} $ denote the value of $\alpha$ at
point $p$, as well as the values of $\alpha$ in element $e$ and the
volume of element $e$, and the sum extends over all the elements
surrounding point $p$.

\subsection{Simple Point/Element/Point Averaging}

In this case, the values of $\alpha$ are cycled between elements 
and points. When going from point values to element values, 
a simple average is taken:
\begin{equation} \label{eq:simple_averaging}
    \alpha_{e} = { 1 \over n_{e} } \sum_i \alpha_i
\end{equation}
where $n_{e}$ denotes the number of nodes (degrees of freedom) of
an element and the sum extends over all the nodes of the element.
After obtaining the new element values via equation \eqref{eq:simple_averaging} the point
averages are again evaluated via equation \eqref{eq:averaging}.
This form of averaging is very crude, but works surprisingly well.

\subsection{$H^1$ (Weak) Laplacian Smoothing}

In this case, the initial values $\alpha_0$ obtained for $\alpha$ 
are smoothed via:
\begin{equation} \label{eq:weak_laplacian_smooth_1}
    \left[ 1 - \lambda \nabla^2 \right] \alpha = \alpha_0 ~~, \quad
   \left. \alpha_{,n} \right|_{\Gamma} = 0
\end{equation}
Here $\lambda$ is a free parameter which may be problem and mesh
dependent (its dimensional value is length squared). Discretization
via finite elements yields:
\begin{equation} \label{eq:weak_laplacian_smooth_2}
    \left[ \Mmat_c + \lambda \Kmat_d \right] \alphavec = 
          \Mmat_{p1p0} \alphavec_0
\end{equation}
where $\Mmat_c, \Kmat_d, \Mmat_{p1p0}$ denote the consistent mass 
matrix, the stiffness or `diffusion' matrix obtained for the Laplacian 
operator and the projection matrix from element values ($\alphavec_0$)
to point values ($\alphavec$).

\subsection{Pseudo-Laplacian Smoothing}

One can avoid the dimensional dependency of $\lambda$ by smoothing via:
\begin{equation} \label{eq:pseudo_laplacian_smooth_1}
    \left[ 1 - \lambda \nabla h^2 \nabla \right] \alpha = \alpha_0
\end{equation}
where $h$ is a characteristic element size. For linear elements, one
can show that this is equivalent to:
\begin{equation} \label{eq:pseudo_laplacian_smooth_2}
    \left[ \Mmat_c + \lambda \left(\Mmat_l - \Mmat_c \right) \right] 
      \alphavec = \Mmat_{p1p0} \alphavec_0
\end{equation}
where $\Mmat_l$ denotes the lumped mass matrix \cite{lohner2008applied}. In the examples shown below this form of smoothing was used for the gradients, setting $\lambda=0.05$.

\section{Examples} \label{sec:examples}

All the numerical examples were carried out using two finite
element codes. The first, FEELAST \cite{lohner2023feelast}, is a 
finite element code based on simple linear (truss), triangular
(plate) and tetrahedral (volume) elements with constant material
properties per element that only solves the linear elasticity
equations. The second, CALCULIX \cite{dhondt2022calculix}, is a 
general, open source
finite element code for structural mechanical applications with 
many element types, material models and options.
The optimization loops were steered via a simple shell-script for 
the adjoint-based optimization. In all cases, a `target' distribution 
of $\alpha(\xvec)$ was given, together with defined external forces 
$\fvec_{\Gamma}$.
The problem was then solved, i.e. the deformations $\uvec(\xvec)$
and strains $\svec(\xvec)$ were obtained and recorded at the
`measurement locations' $\xvec_j,~j=1,m$. This then yielded 
the `measurement pair' $\fvec, \uvec_j,~j=1,m$ or
$\fvec, \svec_j,~j=1,m$ that was used to determine the material
strength distributions $\alpha(\xvec)$ in the field.

The first cases serve to verify that the procedure can recover a uniform strength factor, starting for an arbitrary distribution. The subsequent cases treat the more realistic scenario of trying to determine regions of weakening materials.

\subsection{Crane}

\begin{figure}[!hbt]
    \centering
    \includegraphics[width=0.7\textwidth]{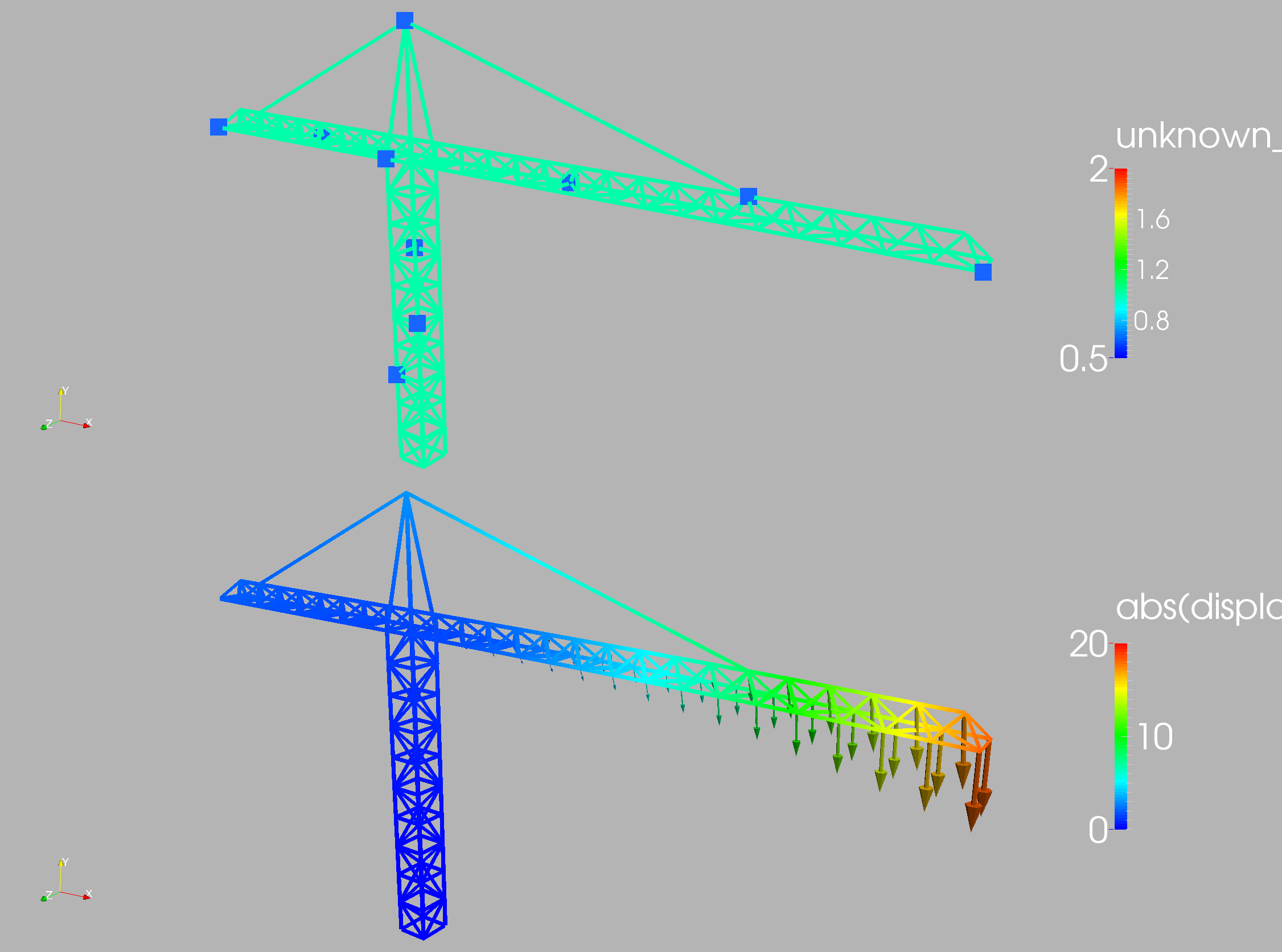}
    \caption{Crane: base case ($\alpha=1.0$)}
    \label{fig:crane_base}
\end{figure}

The case is shown in Figure~\ref{fig:crane_base} and considers a typical crane used
at construction sites. The crane has a height of 1,400~cm, and the
arm has a length of 2,500~cm. A typical truss is about 100~cm long
and has an area of 5~sqcm.
Density, Young's modulus and Poisson rate were set to
$\rho=7.8, E=2 \cdot 10^{12}, \nu=0.3$ respectively (all cgs units).
The two end points on the arm had loads of $f_y=-2.0 \cdot 10^9~gr~cm/sec^2$ 
applied, while the two end points on balancing/back part of the arm
had loads of $f_y=-1.0 \cdot 10^9~gr~cm/sec^2$.
The finite element discretization consisted
of 350 linear truss elements. The loads lead to a deformation
in the vertical direction $w_y=-18~cm$ at the tip of the arm. The top 
figure shows the strength factor $\alpha$ and
the ten measuring points used (which in this case coincide with nodes 
of the finite element mesh), while the bottom figure displays the 
deformation field.

\subsubsection{Displacement Measurements}

\begin{figure}[!hbt]
    \centering
    \includegraphics[width=0.7\textwidth]{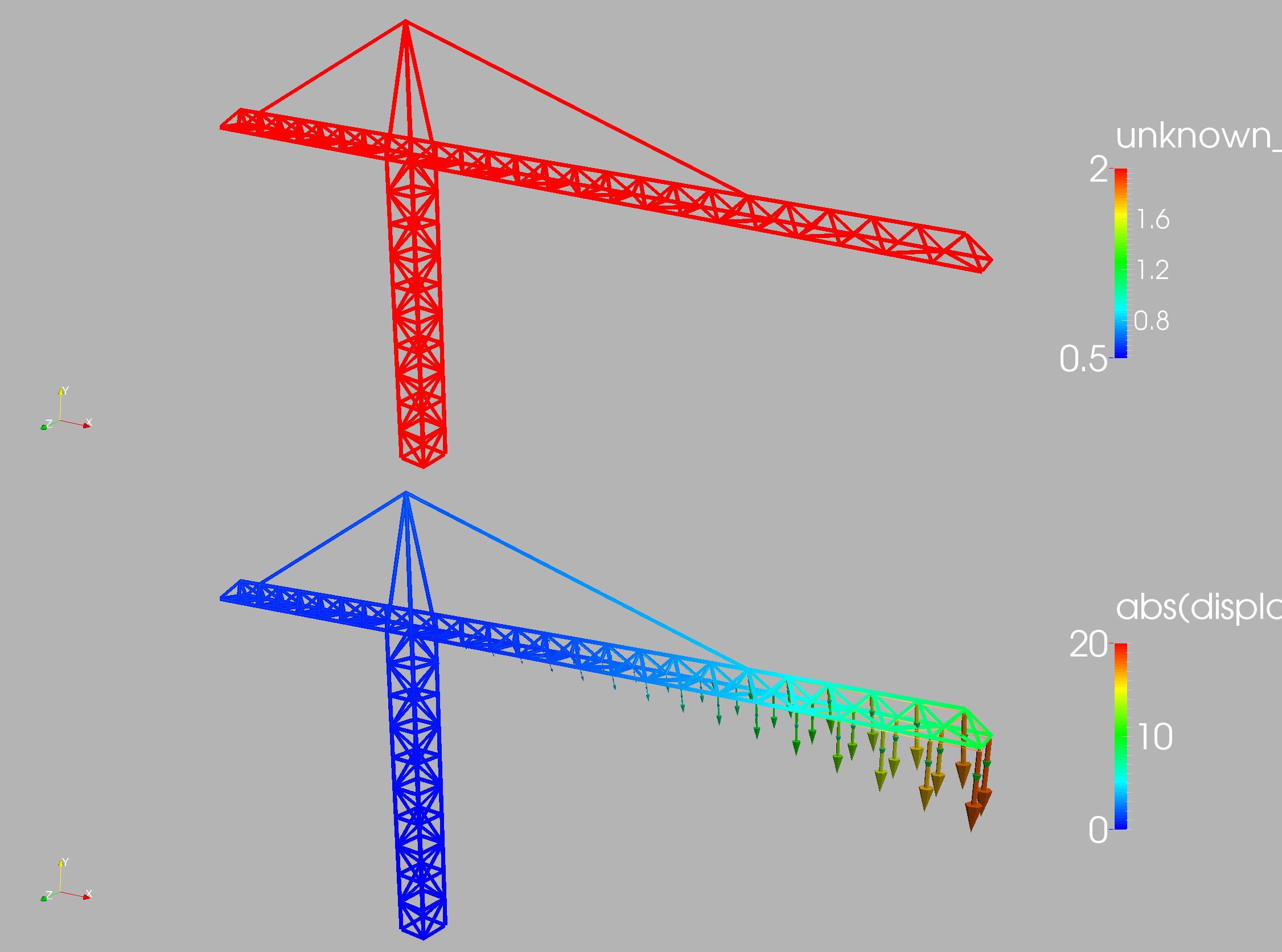}
    \caption{Crane: Start: $\alpha=2.0$, Iteration: 0}
    \label{fig:crane_2a}
\end{figure}

\begin{figure}[!hbt]
    \centering
    \includegraphics[width=0.7\textwidth]{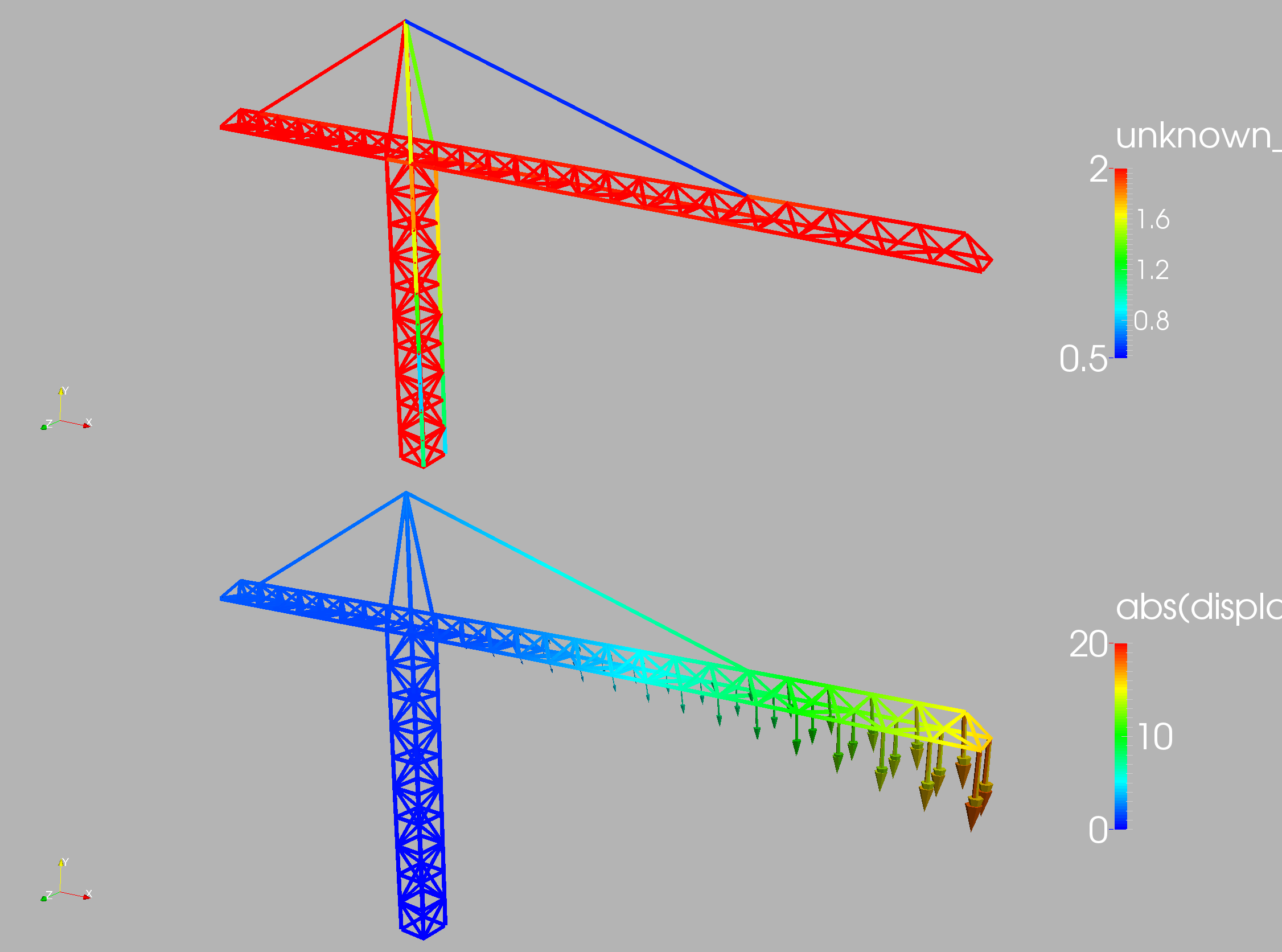}
    \caption{Crane: Start: $\alpha=2.0$, Iteration: 90, No Smoothing
of Gradients}
    \label{fig:crane_2b}
\end{figure}

\begin{figure}[!hbt]
    \centering
    \includegraphics[width=0.7\textwidth]{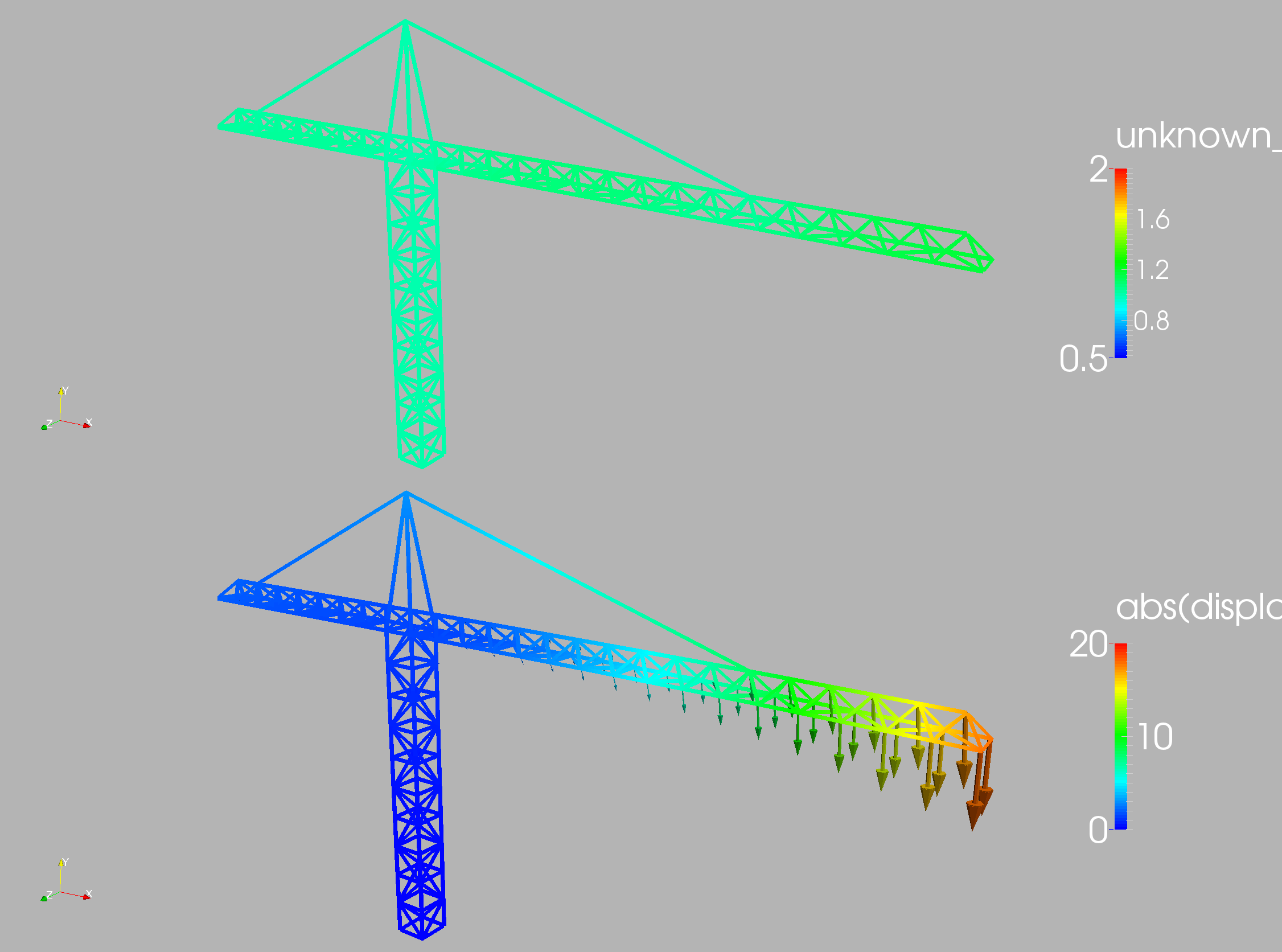}
    \caption{Crane: Start: $\alpha=2.0$, Iteration: 90, Smoothing
of Gradients}
    \label{fig:crane_2c}
\end{figure}

Given the desired/measured displacements at these 10 measuring points,
different starting values for the strength factor $\alpha$ were explored.
Figures~\ref{fig:crane_2a}-\ref{fig:crane_2c} show the results obtained when starting from a uniform
value of $\alpha=2.0$ without Figure~\ref{fig:crane_2b} and with Figure~\ref{fig:crane_2c} gradient smoothing.
One can see that for this case gradient smoothing is essential.

\begin{figure}[!hbt]
    \centering
    \includegraphics[width=0.7\textwidth]{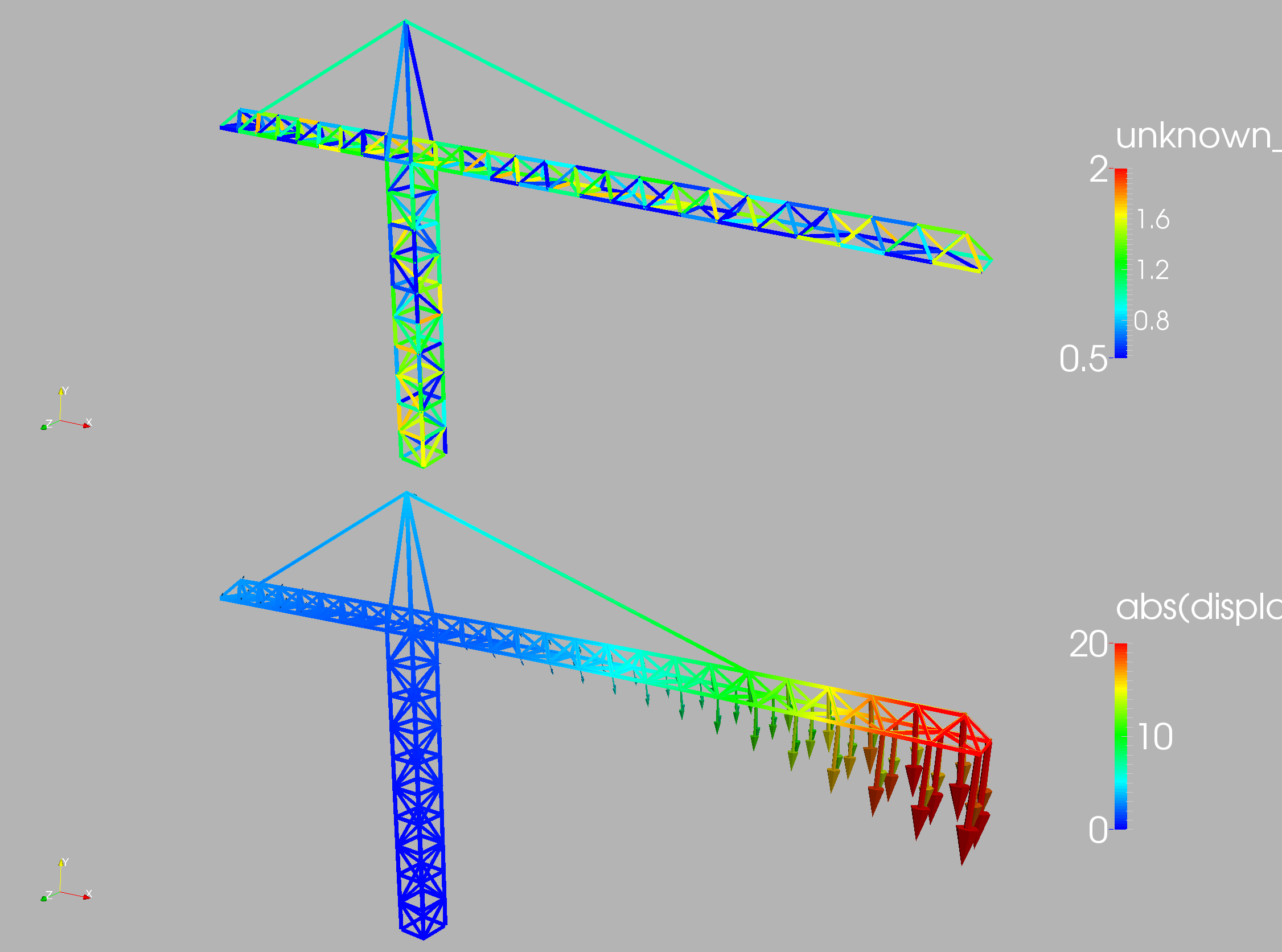}
    \caption{Crane: Start: Random $\alpha$, Iteration: 0}
    \label{fig:crane_3a}
\end{figure}

\begin{figure}[!hbt]
    \centering
    \includegraphics[width=0.7\textwidth]{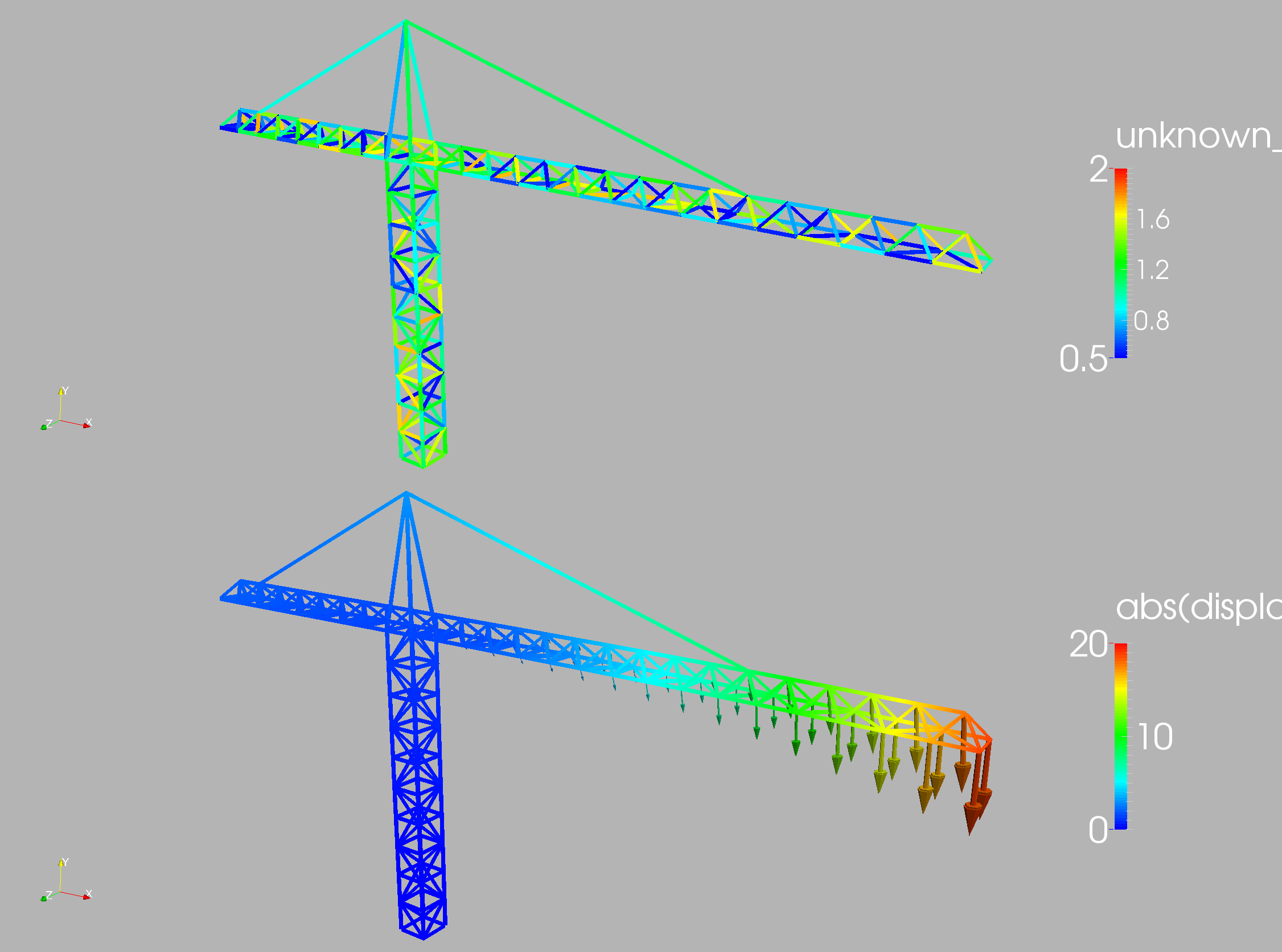}
    \caption{Crane: Start: Random $\alpha$, Iteration: 90, No Smoothing
of Gradients}
    \label{fig:crane_3b}
\end{figure}

\begin{figure}[!hbt]
    \centering
    \includegraphics[width=0.7\textwidth]{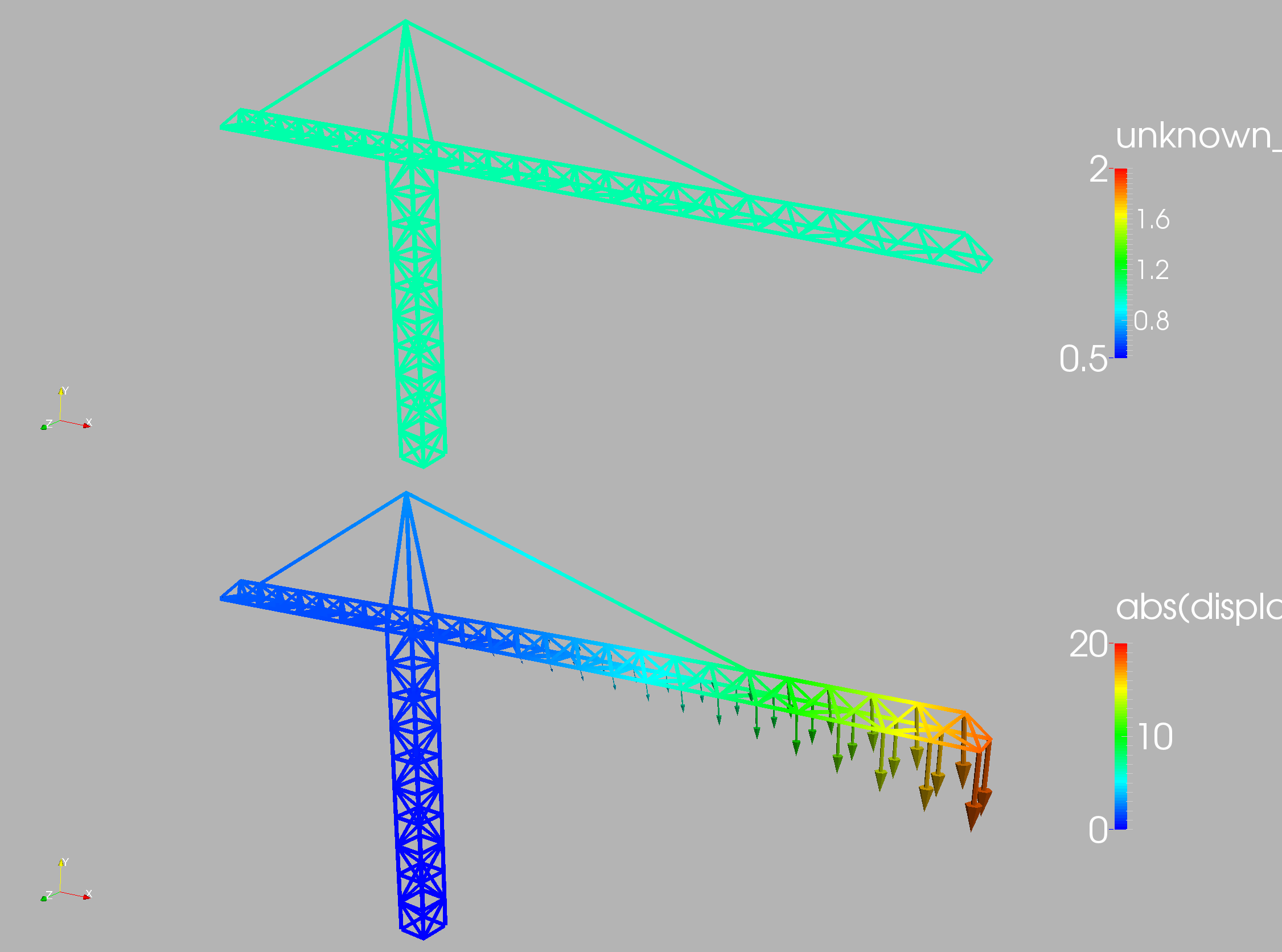}
    \caption{Crane: Start: Random $\alpha$, Iteration: 90, Smoothing
of Gradients}
    \label{fig:crane_3c}
\end{figure}

Figures \ref{fig:crane_3a}, \ref{fig:crane_3b} and \ref{fig:crane_3c} show the results obtained when starting from a random
distribution of $\alpha$ without and with gradient smoothing.
As before, one can see that for this case gradient smoothing is 
essential.

\subsubsection{Strain Measurements}

\begin{figure}[!hbt]
    \centering
    \includegraphics[width=0.7\textwidth]{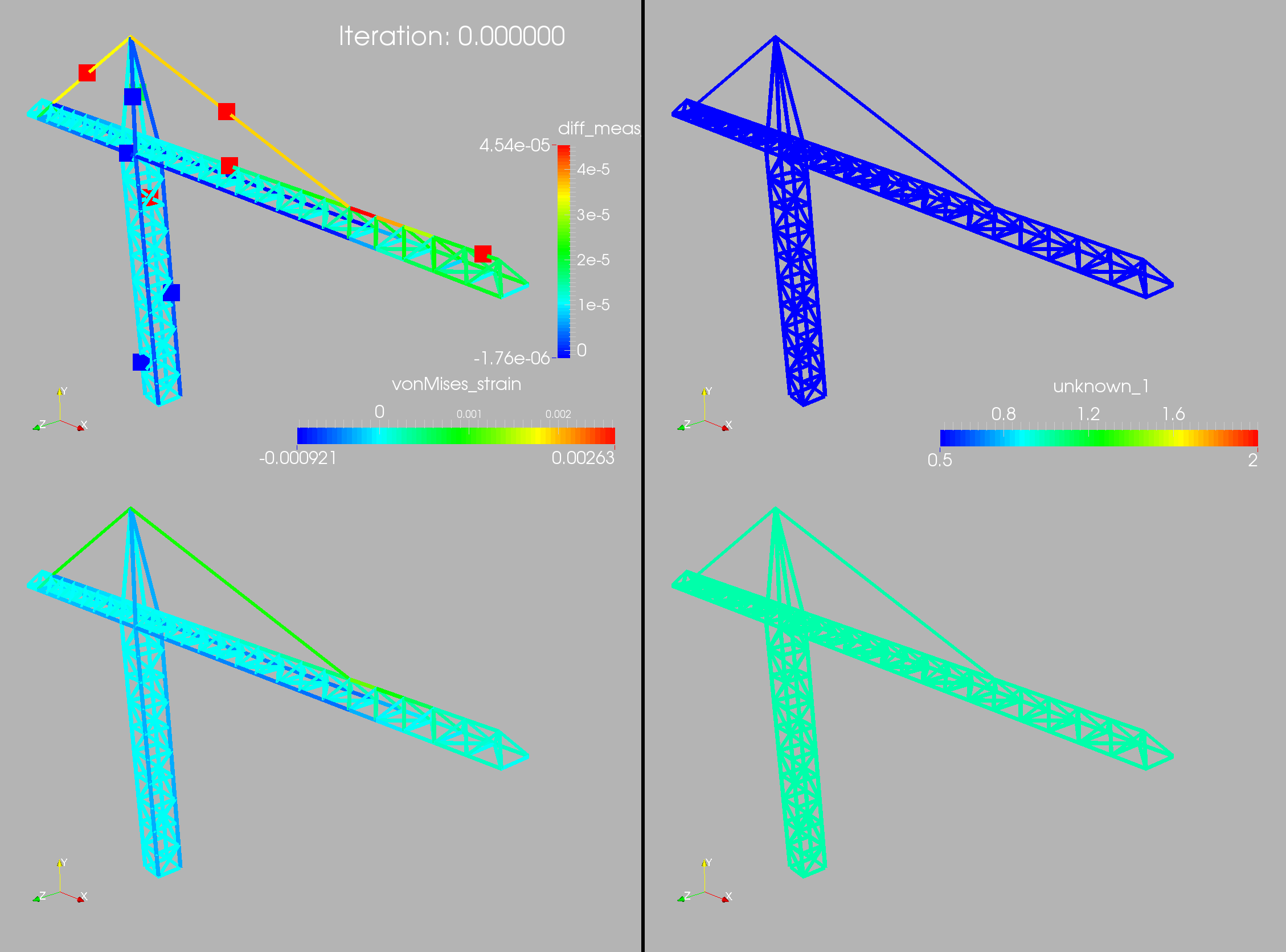}
    \caption{Crane: Start: $\alpha=0.5$, Iteration: 0}
    \label{fig:crane_4a}
\end{figure}

\begin{figure}[!hbt]
    \centering
    \includegraphics[width=0.7\textwidth]{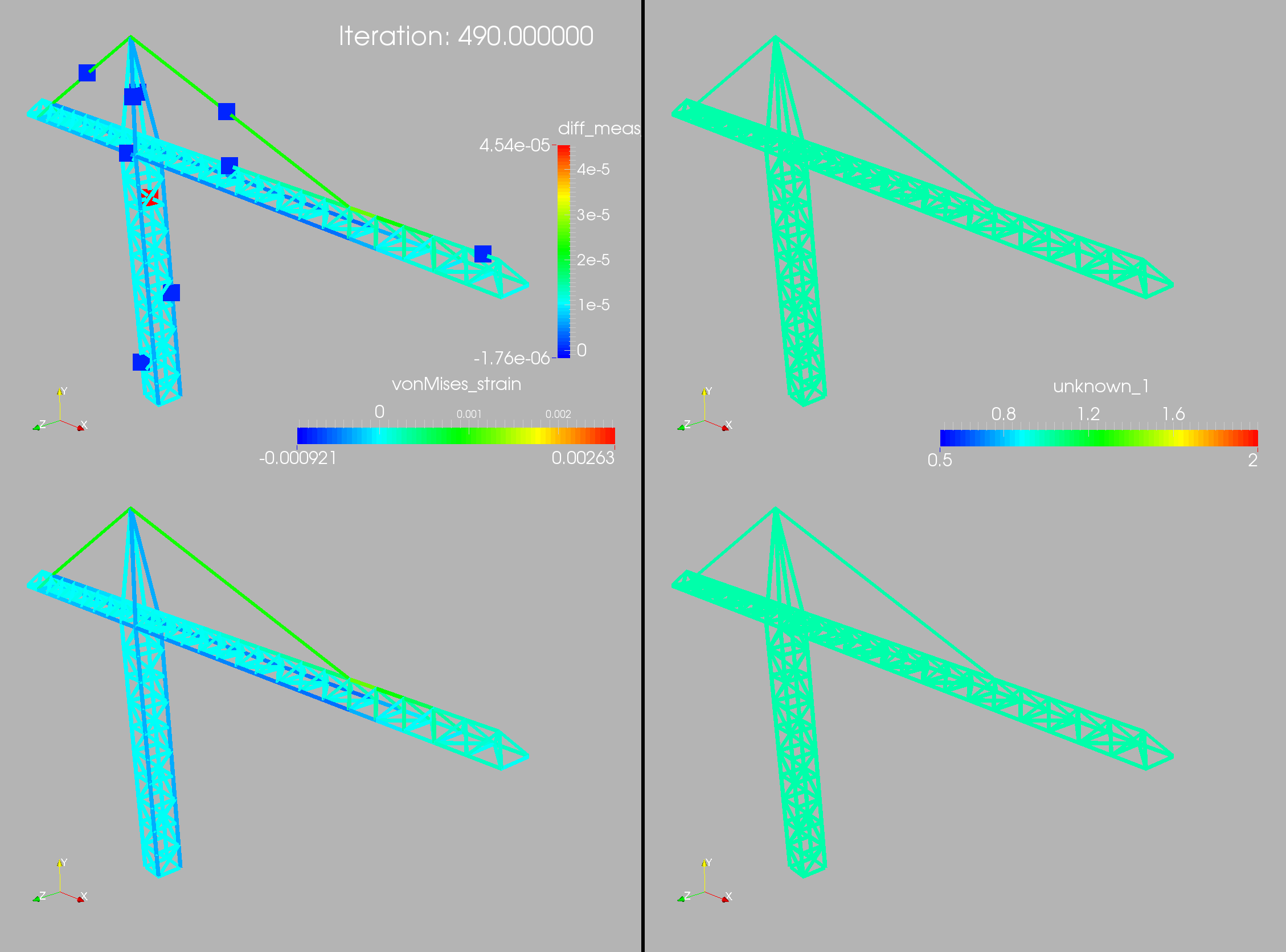}
    \caption{Crane: Start: $\alpha=0.5$, Iteration: 490}
    \label{fig:crane_4b}
\end{figure}

\begin{figure}[!hbt]
    \centering
    \includegraphics[width=0.7\textwidth]{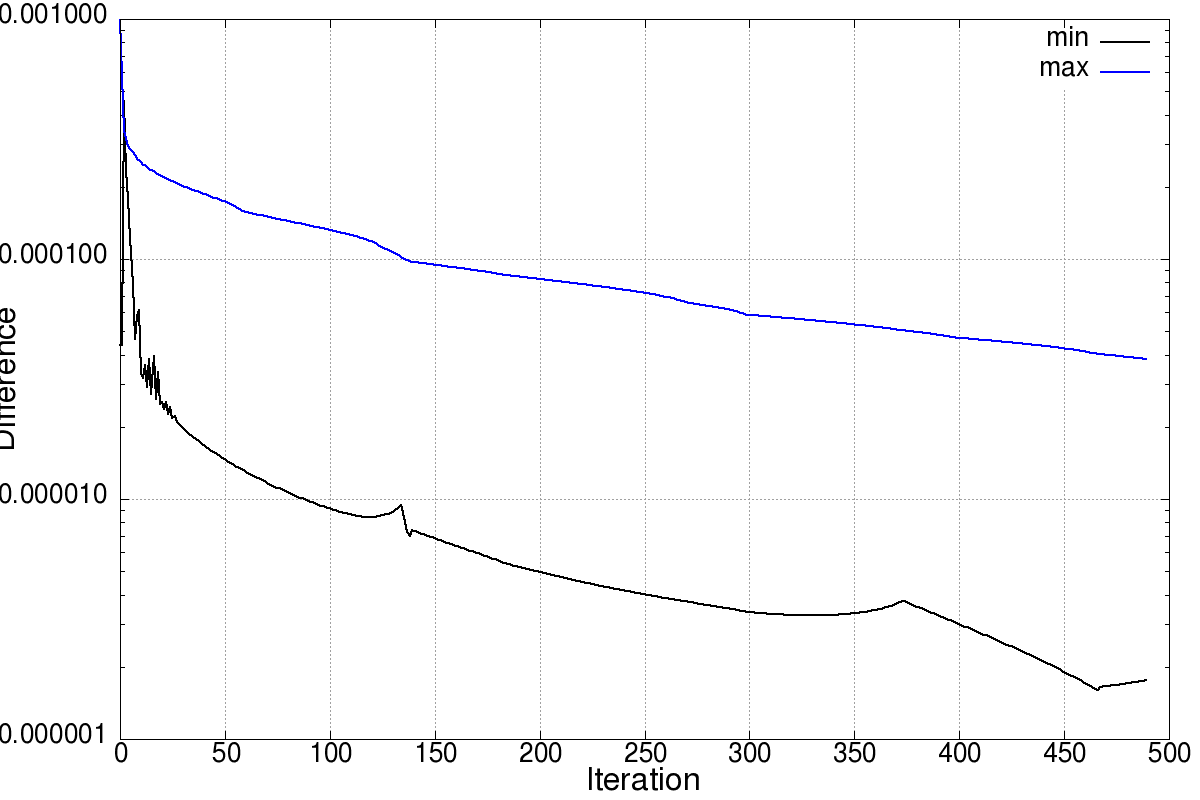}
    \caption{Crane: Start: $\alpha=0.5$, abs(Minimum) and Maximum}
    \label{fig:crane_4c}
\end{figure}

Ten strain measuring points were defined in trusses along the
structure (see top left of Figures~\ref{fig:crane_4a}, \ref{fig:crane_4b} and \ref{fig:crane_4c}).
Given the desired/measured strains at these 10 measuring points,
different starting values for the strength factor $\alpha$ were explored.
The results obtained and behaviours observed were very similar to the
cases with displacement measurements: gradient smoothing was essential.
Therefore, gradient smoothing has always been applied for all the 
results shown in the sequel.
Figures \ref{fig:crane_4a}, \ref{fig:crane_4b} and \ref{fig:crane_4c} show the results obtained when starting from a uniform
value of $\alpha=0.5$. The top figures show the actual values while the
bottom part shows the expected strain and strength distribution in the
trusses. Note also on the top left the differences in target and
actual strain at the measurement points.

\begin{figure}[!hbt]
    \centering
    \includegraphics[width=0.7\textwidth]{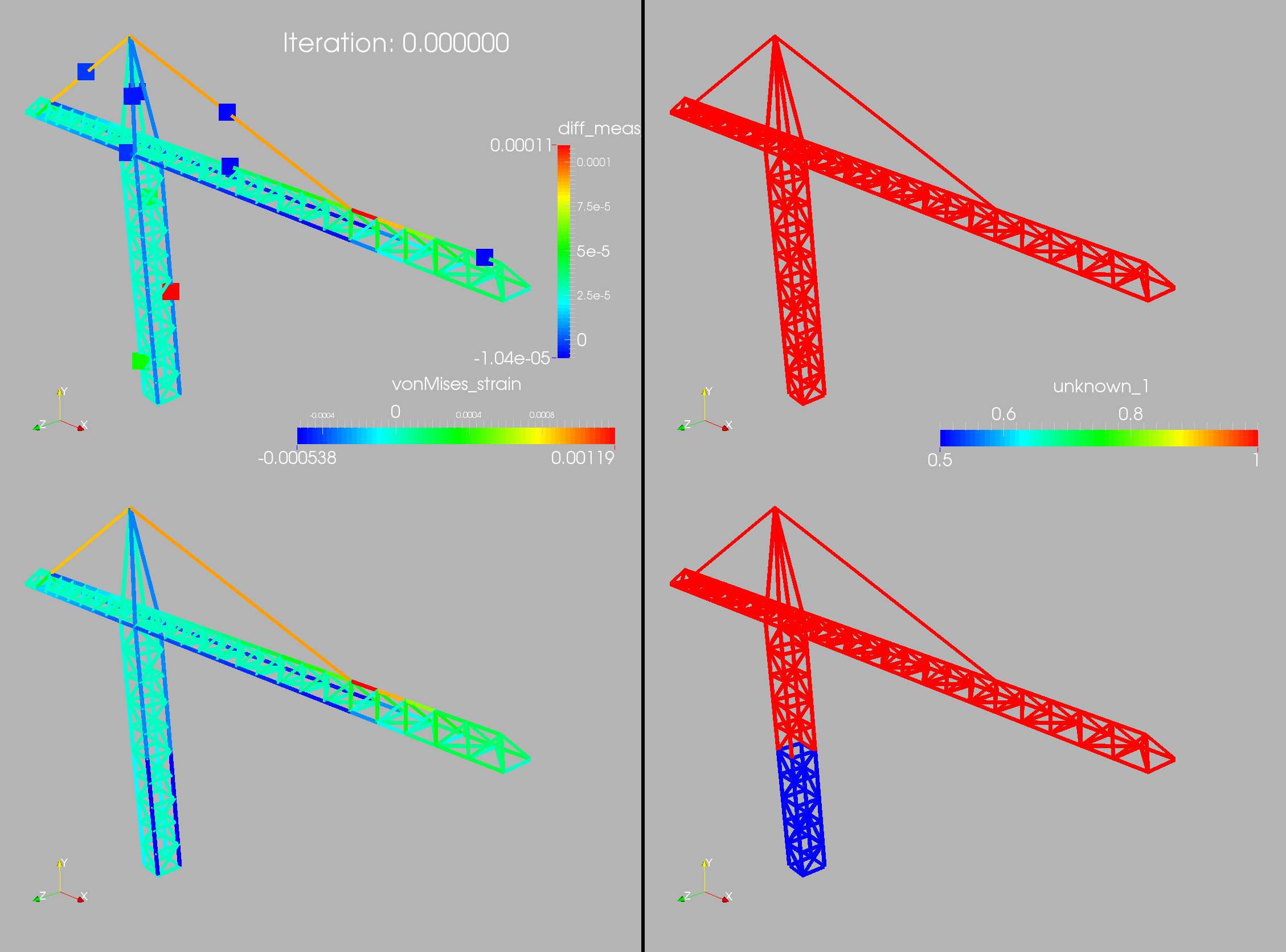}
    \caption{Weakened Crane: Start: $\alpha=1.0$, Iteration: 0}
    \label{fig:crane_5a}
\end{figure}

\begin{figure}[!hbt]
    \centering
    \includegraphics[width=0.7\textwidth]{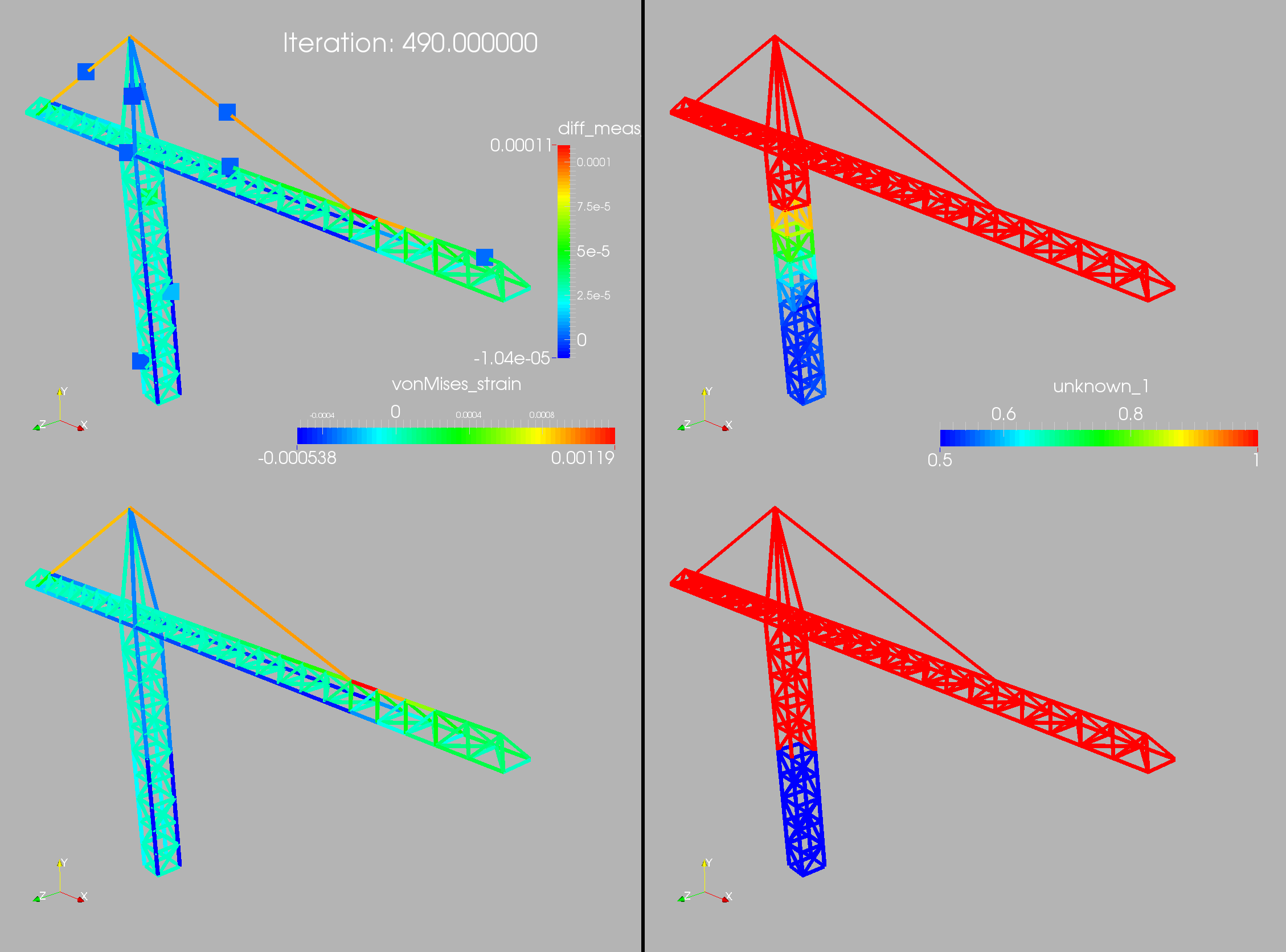}
    \caption{Weakened Crane: Start: $\alpha=1.0$, Iteration: 490}
    \label{fig:crane_5b}
\end{figure}

\begin{figure}[!hbt]
    \centering
    \includegraphics[width=0.7\textwidth]{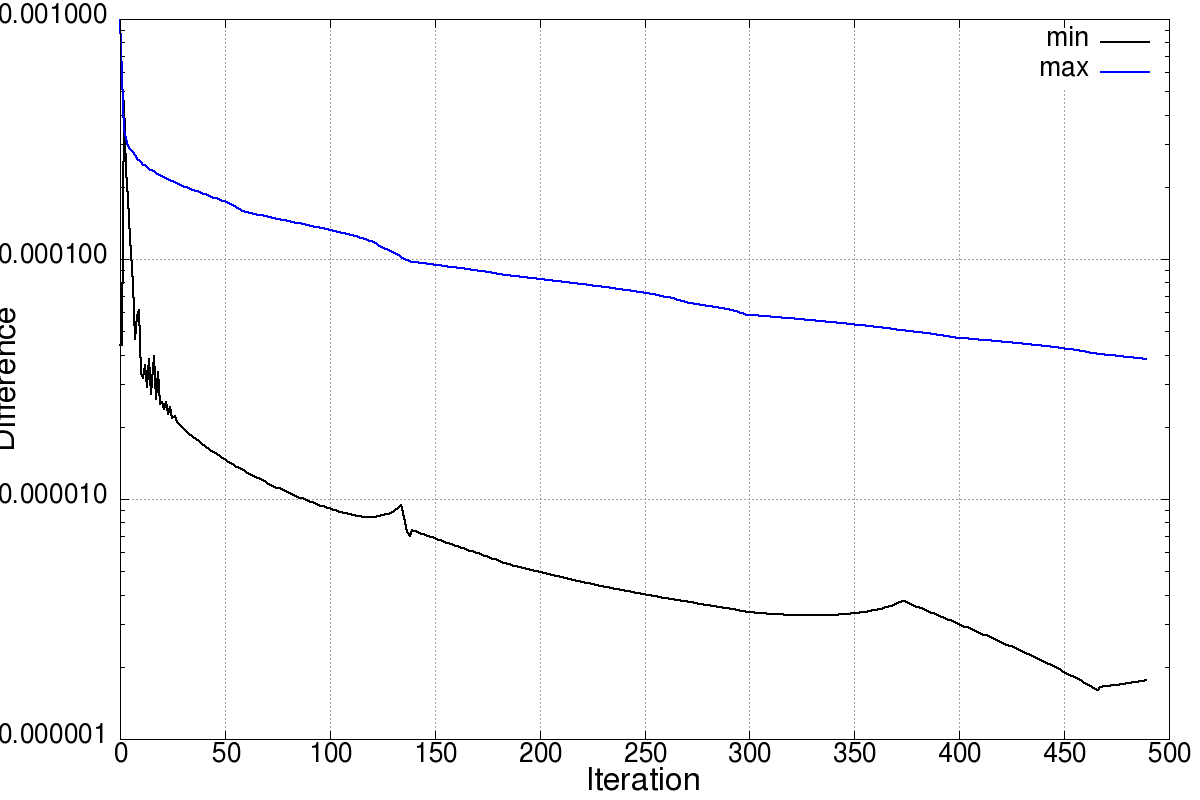}
    \caption{Weakened Crane: Start: $\alpha=1.0$, abs(Minimum) and 
Maximum}
    \label{fig:crane_5c}
\end{figure}

Figures \ref{fig:crane_5a}, \ref{fig:crane_5b} and \ref{fig:crane_5c} show the results obtained when starting from a uniform
value of $\alpha=1.0$ for the case that the lower part of the crane
tower has been weakened to $\alpha=0.5$. As before, the top figures 
show the actual values while the bottom part shows the expected 
strain and strength distribution in the trusses.

\subsubsection{Displacement Measurements With Multiple Loads}

\begin{figure}[!hbt]
    \centering
    \includegraphics[width=0.7\textwidth]{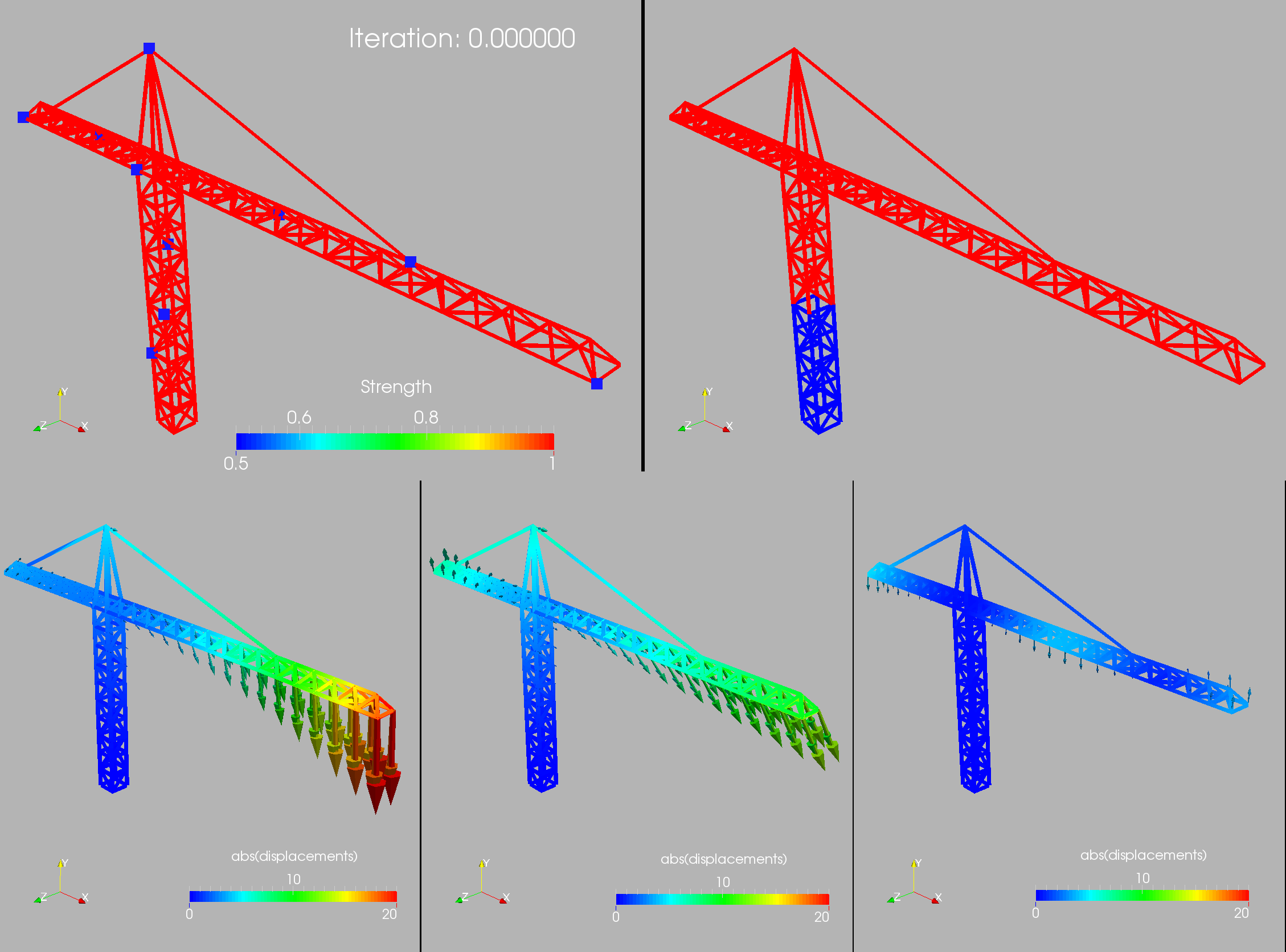}
    \caption{Weakened Crane: Start: $\alpha=1.0$, Iteration: 0}
    \label{fig:crane_6a}
\end{figure}

\begin{figure}[!hbt]
    \centering
    \includegraphics[width=0.7\textwidth]{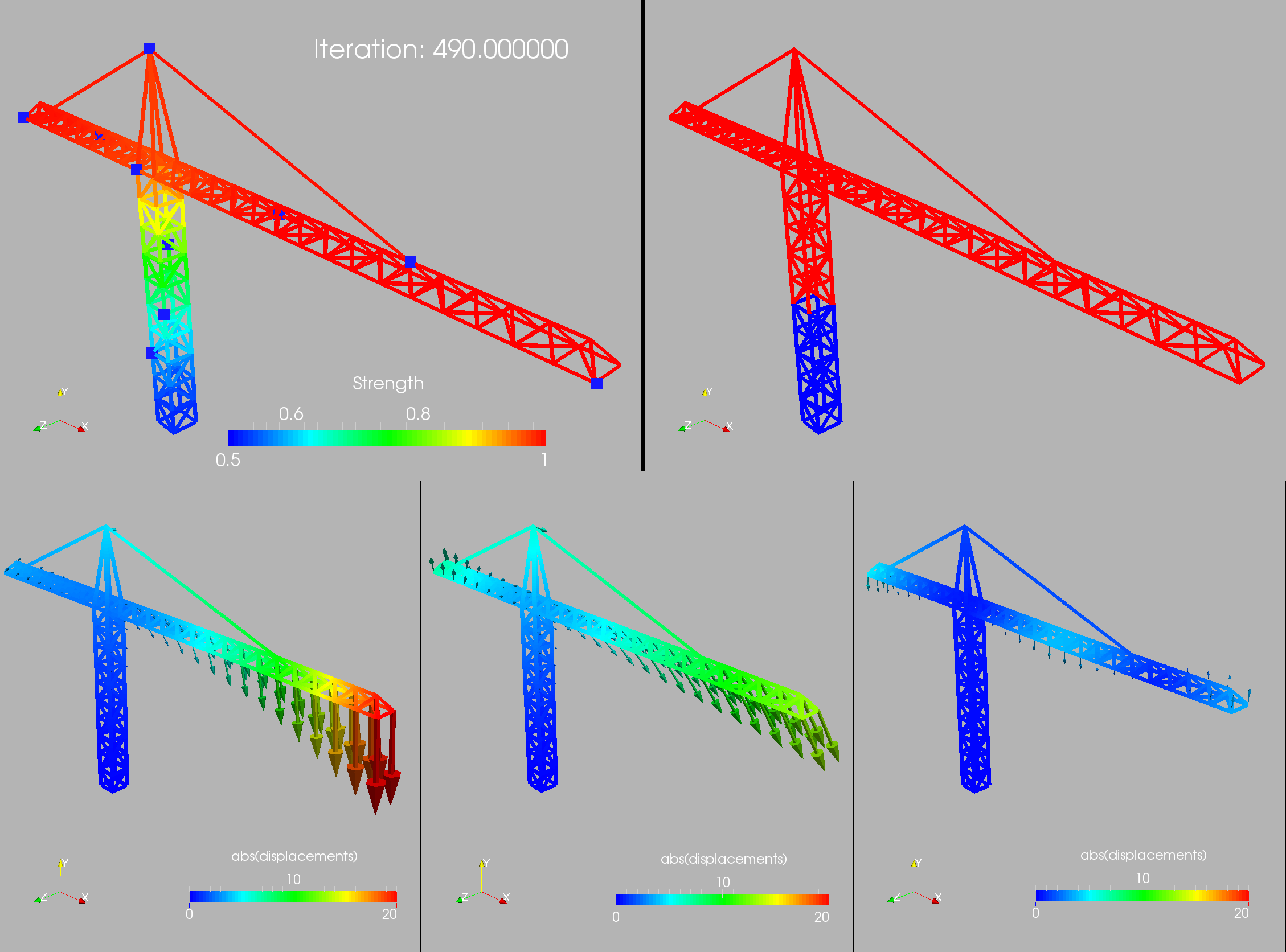}
    \caption{Weakened Crane: Start: $\alpha=1.0$, Iteration: 490}
    \label{fig:crane_6b}
\end{figure}

\begin{figure}[!hbt]
    \centering
    \includegraphics[width=0.7\textwidth]{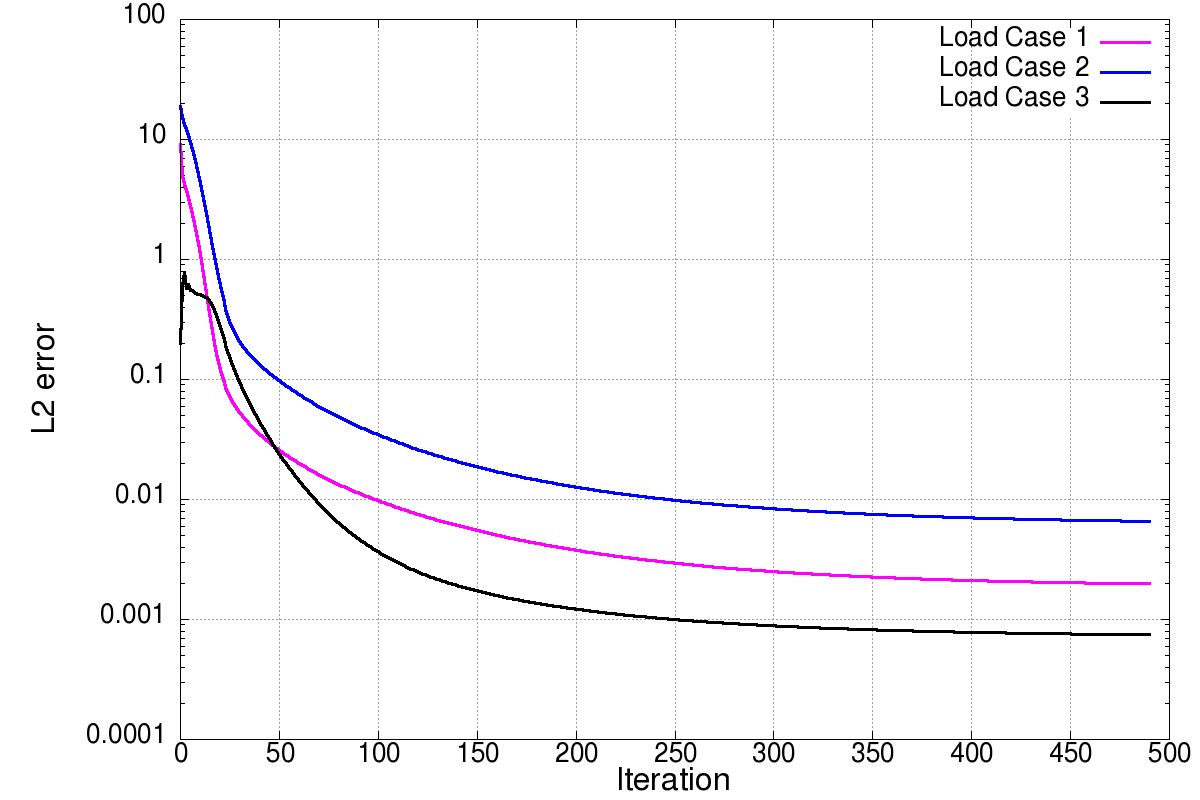}
    \caption{Weakened Crane: Convergence History for the Load Cases}
    \label{fig:crane_6c}
\end{figure}

The same `weakened bottom' scenario was also computed for the 10
displacement measurement points shown before, but with 3 load 
scenarios. The first is the same as before, the second induces a
torsion of the mast and the third applies forces between the
mast and the end of the arm. Figures \ref{fig:crane_6a}, \ref{fig:crane_6b} and \ref{fig:crane_6c} show the results 
obtained when starting from a uniform value of $\alpha=1.0$.
In the figures, the top left shows the computed strength factor,
the top right the desired (exact) strength factor, while the
bottom shows the displacements (computed and desired overlapped) 
for the 3 load cases.

\subsection{Footbridge}

\begin{figure}[!hbt]
    \centering
    \includegraphics[width=0.7\textwidth]{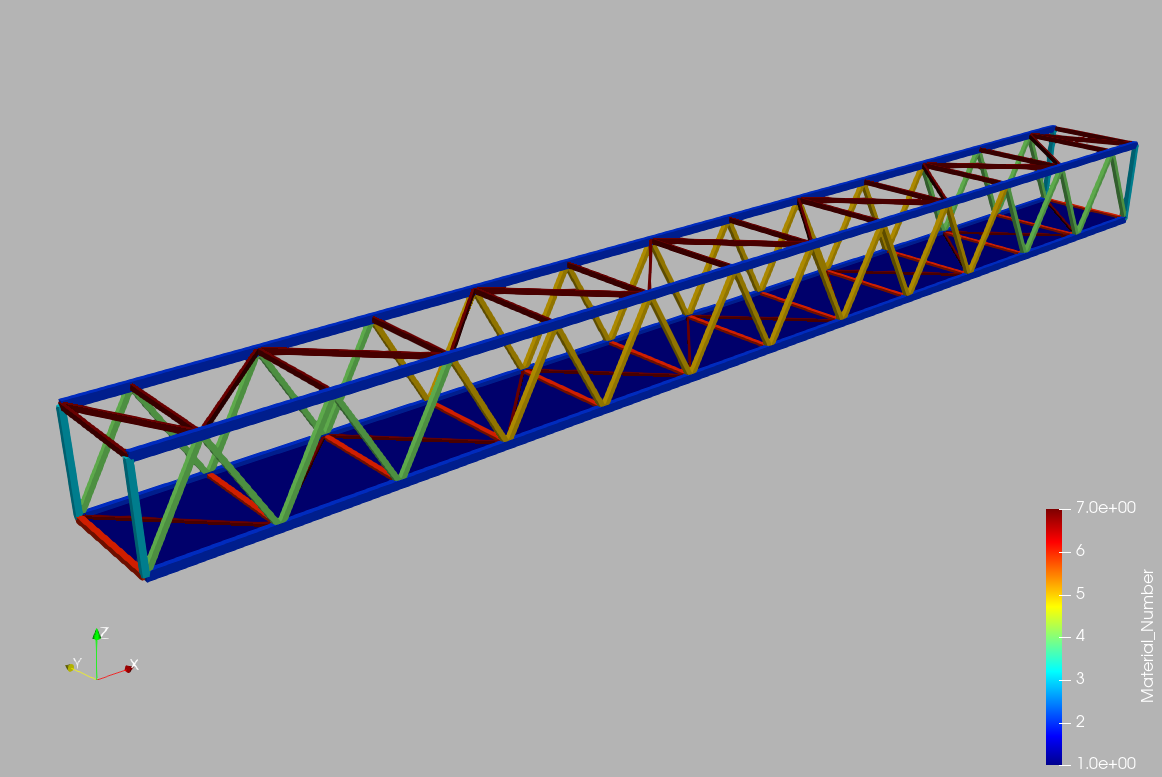}
    \caption{Footbridge: materials}
    \label{fig:footbridge_materials}
\end{figure}

This case considers a typical footbridge and was taken from \cite{kilikevivcius2020influence}.
The forces and material number of the trusses, whose dimensions (all
units in mks) have been compiled in Table \ref{table:footbridge}, can be 
discerned from Figure \ref{fig:footbridge_materials}.
Density, Young's modulus and Poisson rate were set to
$\rho=7,800, E=2 \cdot 10^{11}, \nu=0.3$ respectively.

\begin{table}[!hbt]
\begin{tabular}{@{}ll@{}}
\toprule
\textbf{Material \#} & \textbf{Properties. Dimensions in mm}                  \\ \midrule
1                    & Steel plate. $t = 10$                                 \\
2                    & Steel beam. Hollow section $300 \times 200 \times 12$ \\
3                    & Steel beam. Hollow section $200 \times 200 \times 10$ \\
4                    & Steel beam. Hollow section $180 \times 180 \times 10$ \\
5                    & Steel beam. Hollow section $180 \times 180 \times 5$  \\
6                    & Steel beam. Hollow section $200 \times 200 \times 10$ \\
7                    & Steel beam. Hollow section $200 \times 100 \times 5$  \\ \bottomrule
\end{tabular}
\caption{Materials and dimensions of plates and beams.}
\label{table:footbridge}
\end{table}

\begin{figure}[!hbt]
    \centering
    \subfigure[Target displacements and sensor locations]{\includegraphics[width=0.48\linewidth]{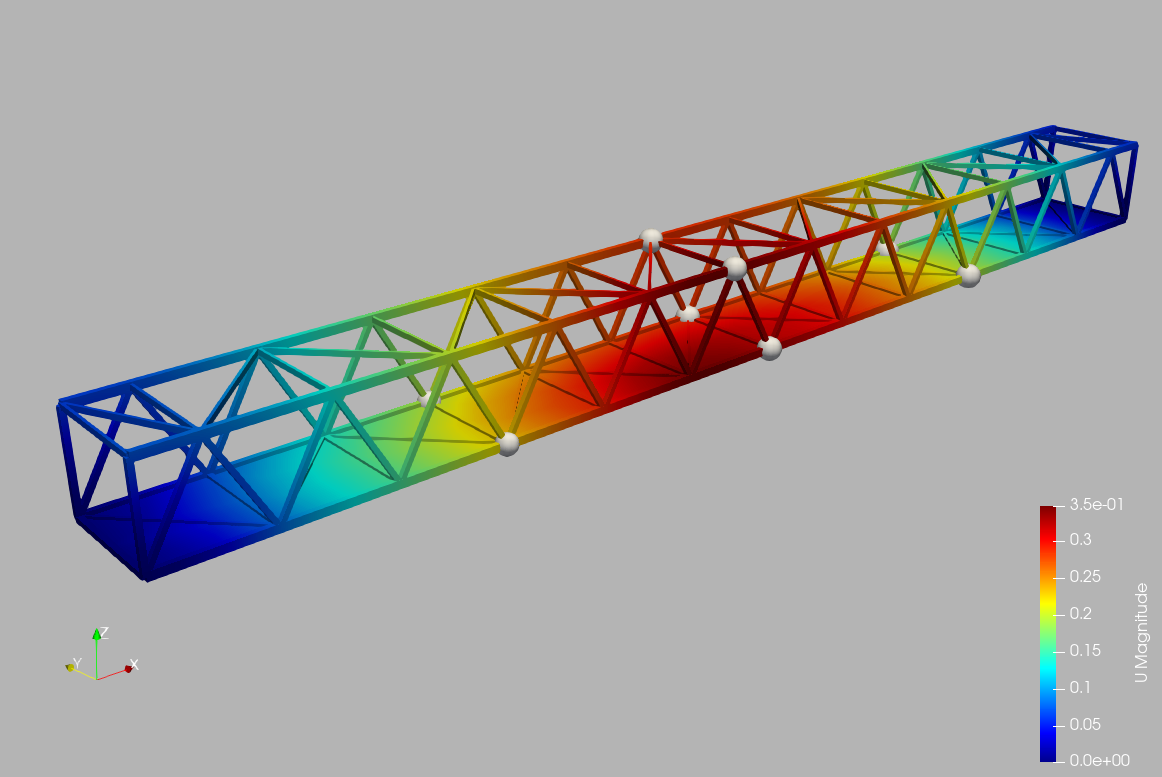}}
    \subfigure[Target strength factor]{\includegraphics[width=0.48\linewidth]{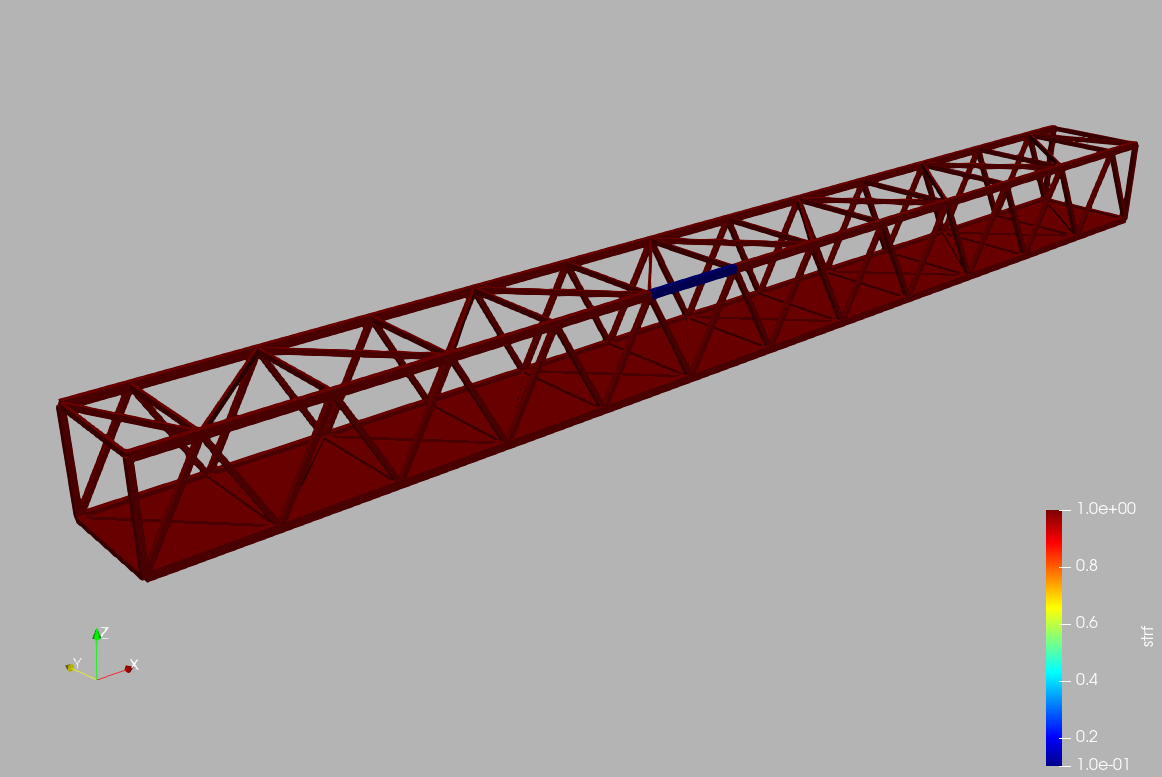}}
    \caption{Footbridge: Target conditions}
    \label{fig:footbridge_target}
\end{figure}

\begin{figure}[!hbt]
    \centering
    \subfigure[Initial displacements and sensor values]{\includegraphics[width=0.48\linewidth]{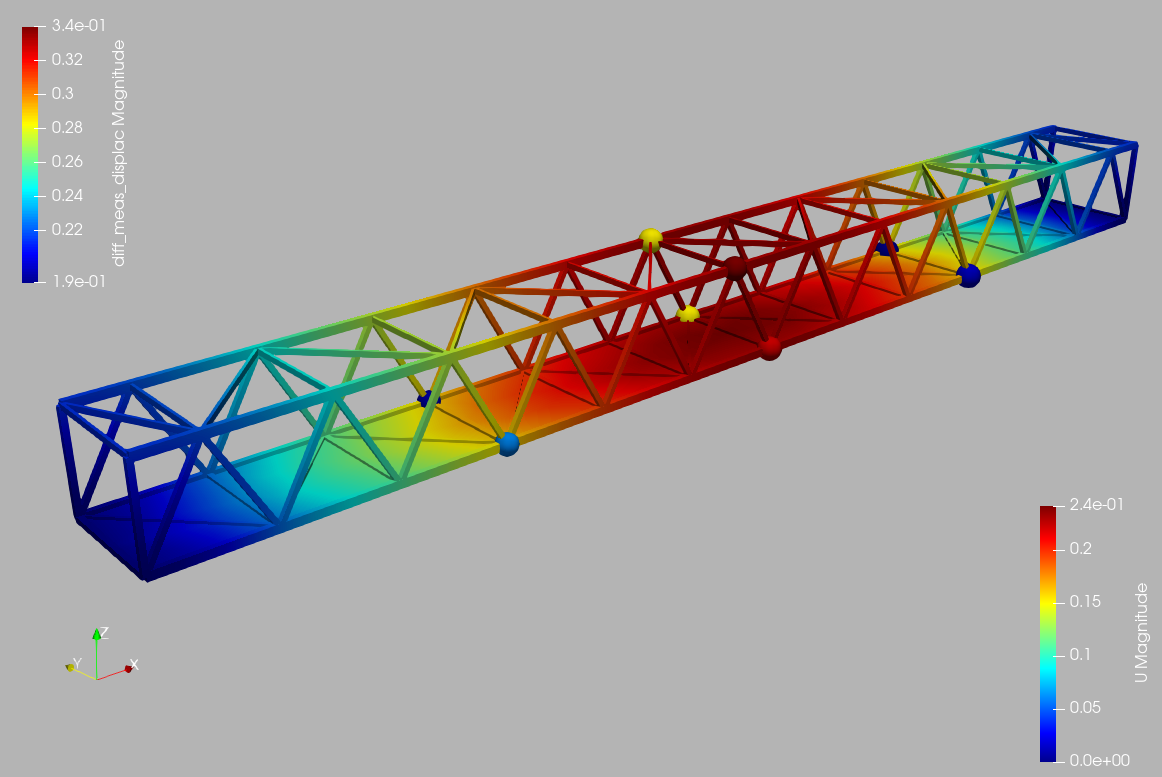}}
    \subfigure[Initial strength factor]{\includegraphics[width=0.48\linewidth]{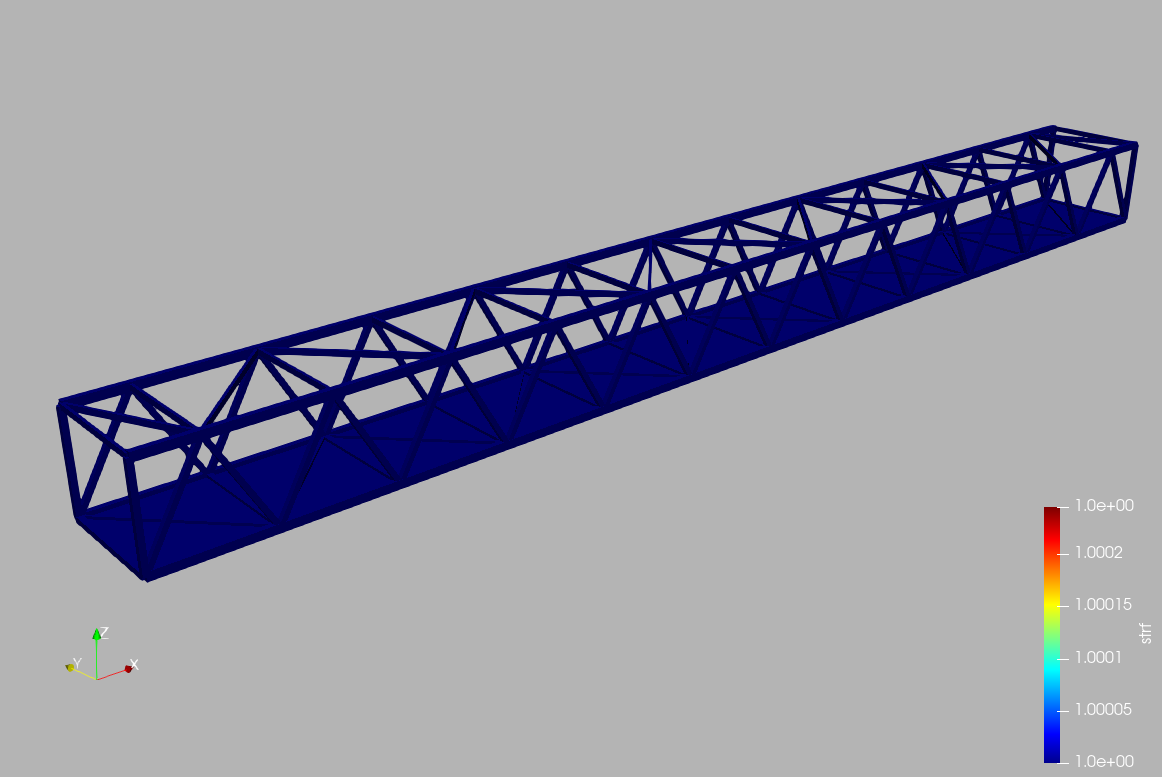}}
    \caption{Footbridge: Initial conditions}
    \label{fig:footbridge_initial}
\end{figure}

\begin{figure}[!hbt]
    \centering
    \subfigure[Displacements obtained]{\includegraphics[width=0.48\linewidth]{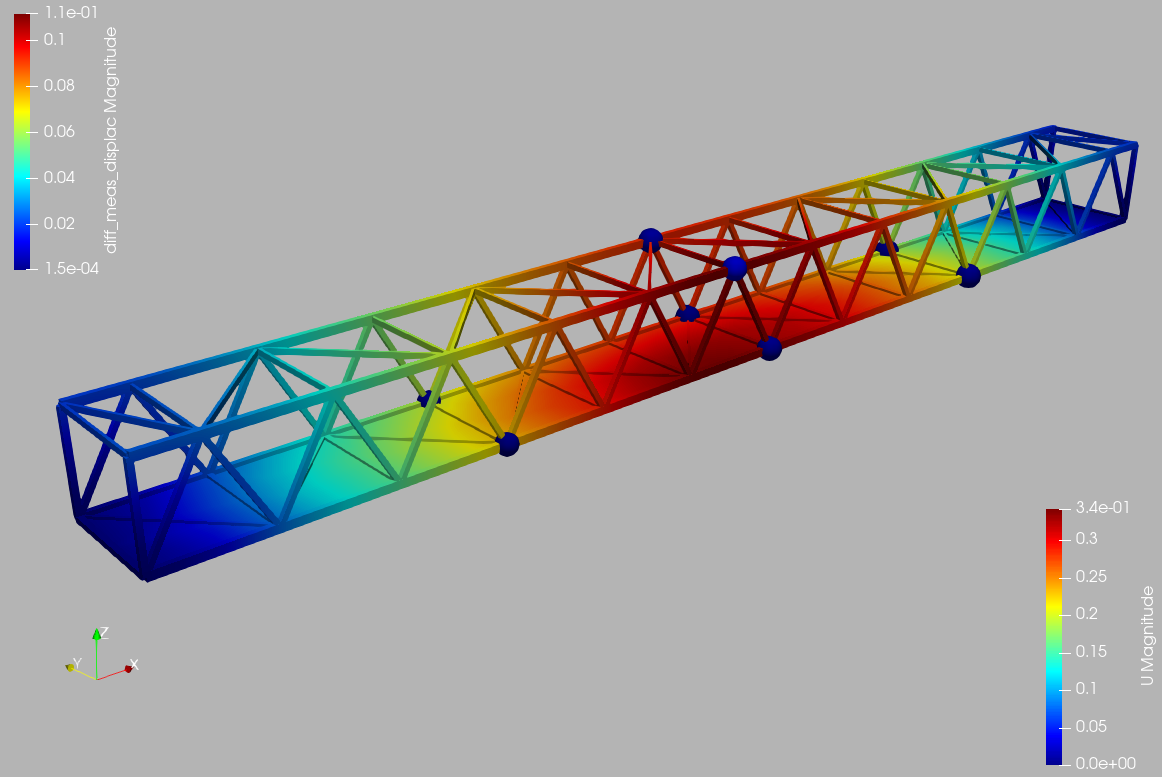}}
    \subfigure[Strength factor obtained]{\includegraphics[width=0.48\linewidth]{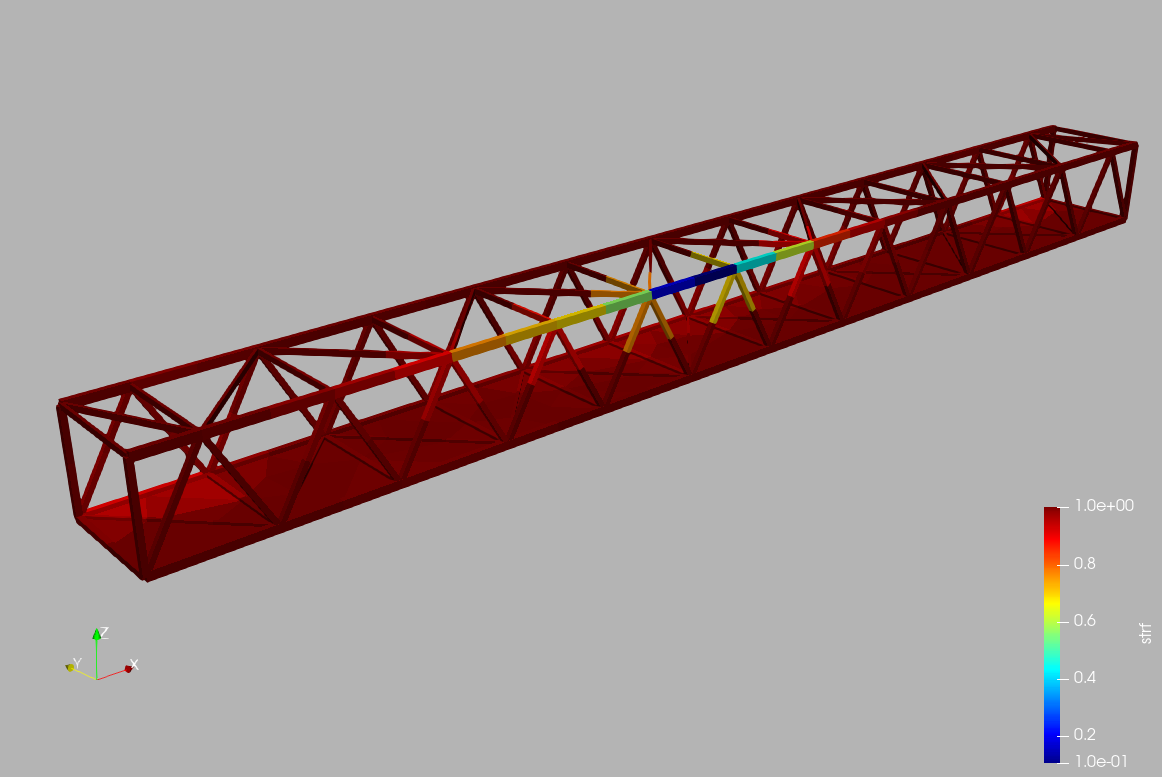}}
    \caption{Footbridge: Solution obtained}
    \label{fig:footbridge_optim}
\end{figure}

\begin{figure}
    \centering
    \includegraphics[width=0.8\textwidth]{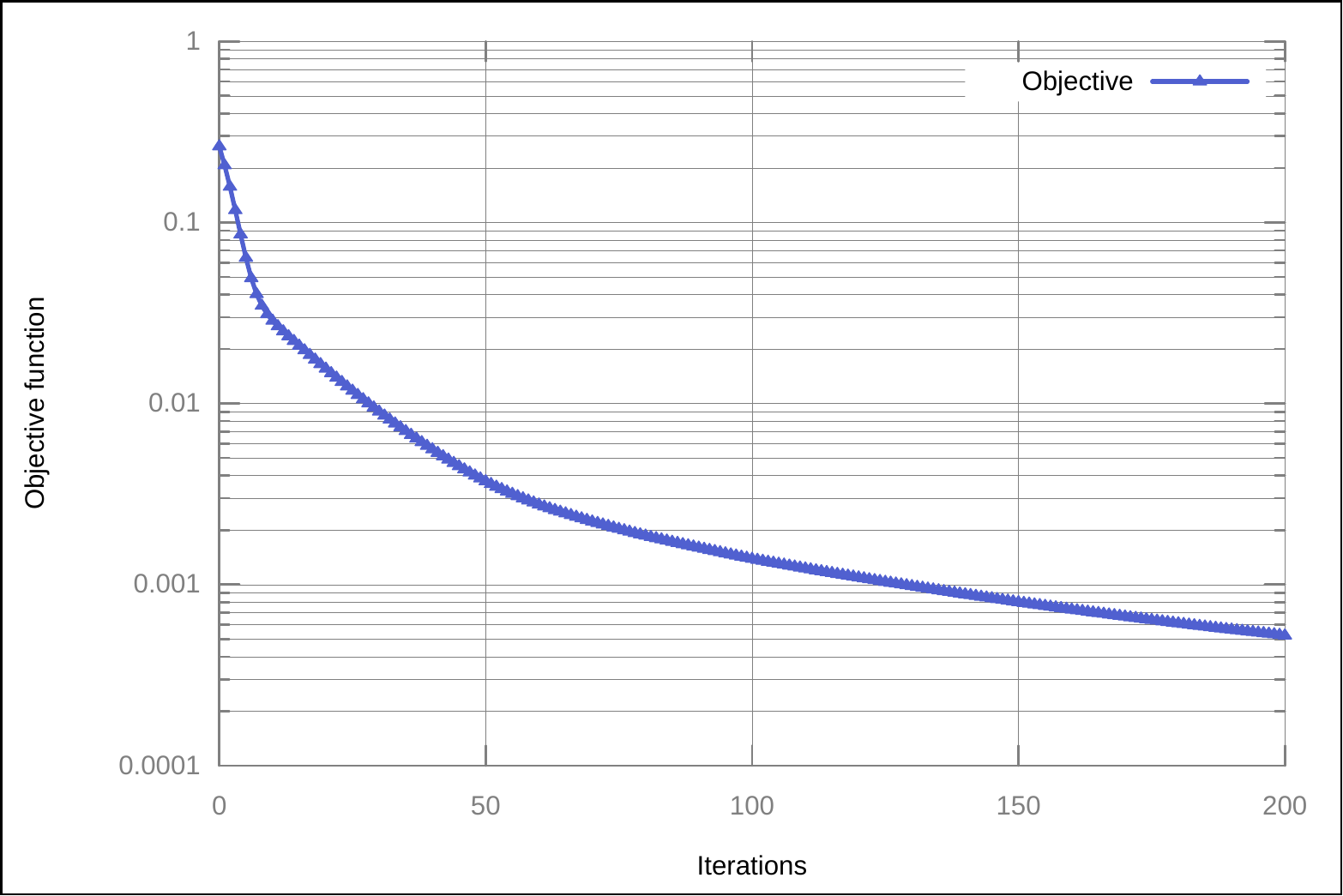}
    \caption{Footbridge: Objective function history}
    \label{fig:footbridge_obj}
\end{figure}

The structure was modeled using 136 shell and 329 beam elements. 
The bridge is under a distributed load of $1$ MPa on the downwards direction, 
applied to every plate, as well as gravity. 
Figure \ref{fig:footbridge_target} shows the target case where we 
have $E=0.1E_0$ at one beams in the structure. 
In \ref{fig:footbridge_target}(a), the location of the 8 sensors is shown, 
along with the target displacements. Starting from a uniform value of 
$\alpha=1.0$, shown in Figure \ref{fig:footbridge_initial}, the
target case is nearly reproduced in 200 steepest descent iterations, 
as can be seen in Figure \ref{fig:footbridge_optim}. The evolution of the 
objective function is shown in figure \ref{fig:footbridge_obj}.

\subsection{Plate with Hole}

\begin{figure}
    \centering
    \includegraphics[width=0.7\textwidth]{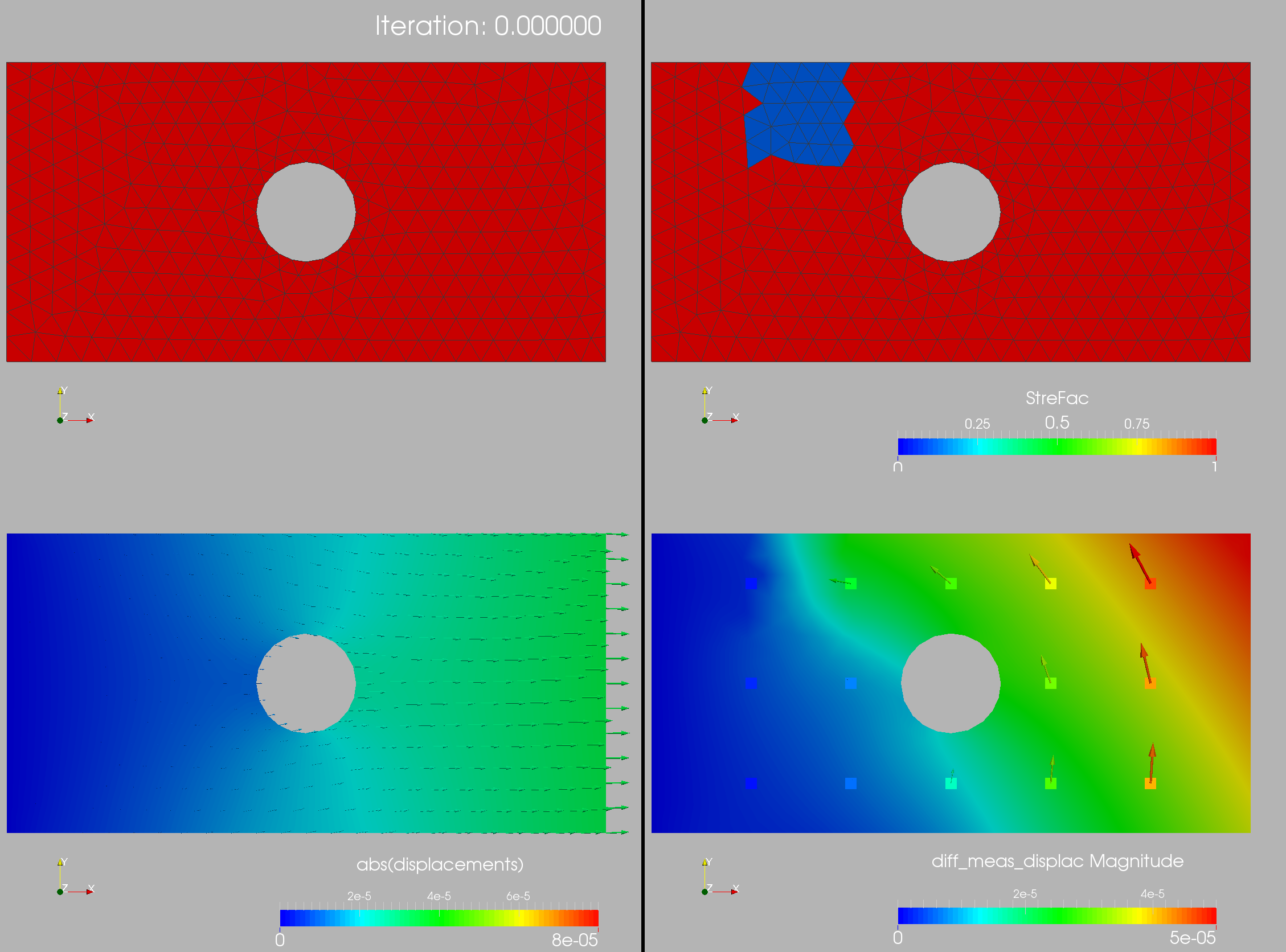}
    \caption{Plate With Hole: Start: $\alpha=1.0$, Iteration: 0}
    \label{fig:plate_with_hole_1}
\end{figure}

\begin{figure}
    \centering
    \includegraphics[width=0.7\textwidth]{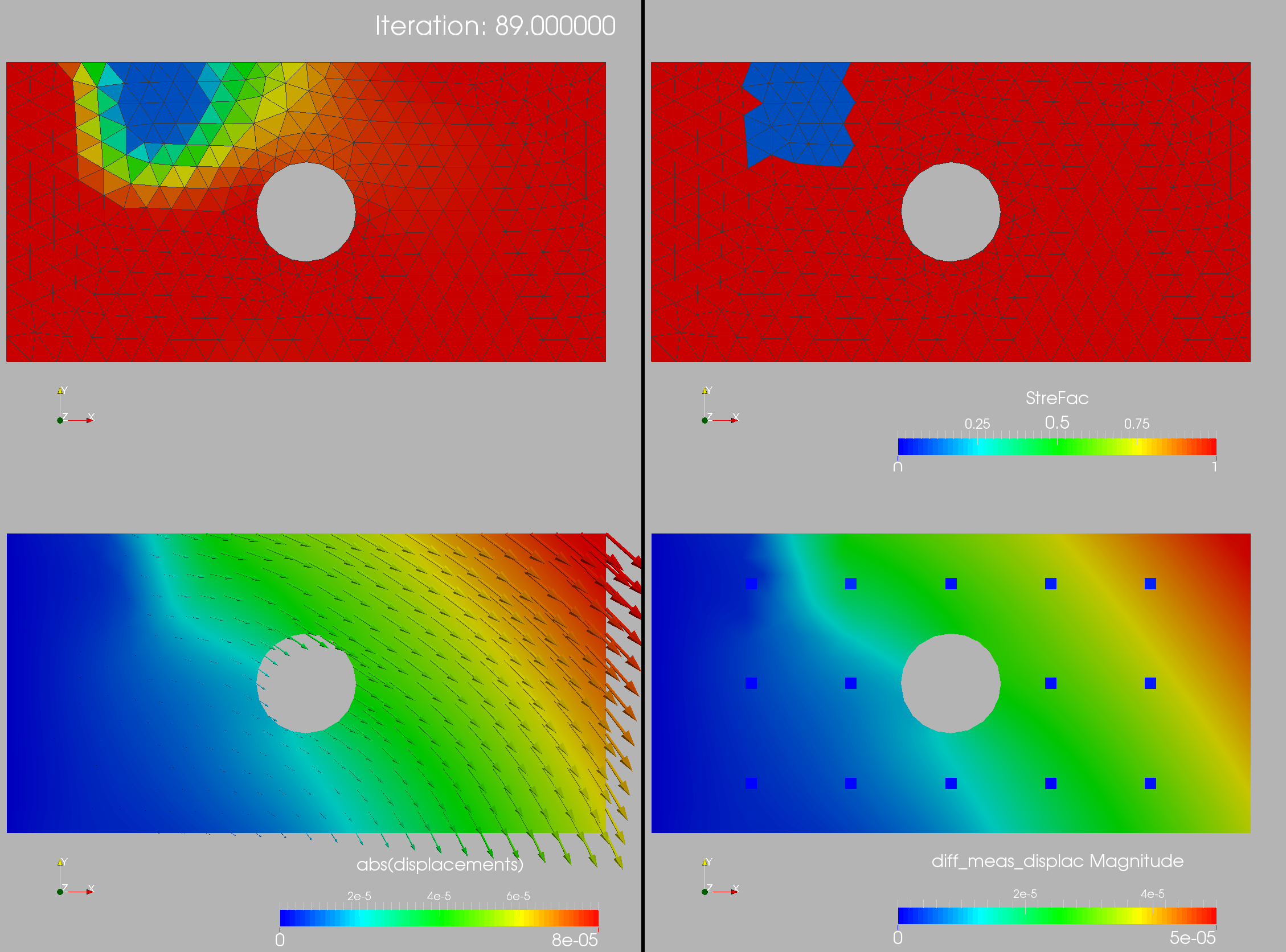}
    \caption{Plate With Hole: Start: $\alpha=1.0$, Iteration: 89}
    \label{fig:plate_with_hole_2}
\end{figure}

\begin{figure}
    \centering
    \includegraphics[width=0.7\textwidth]{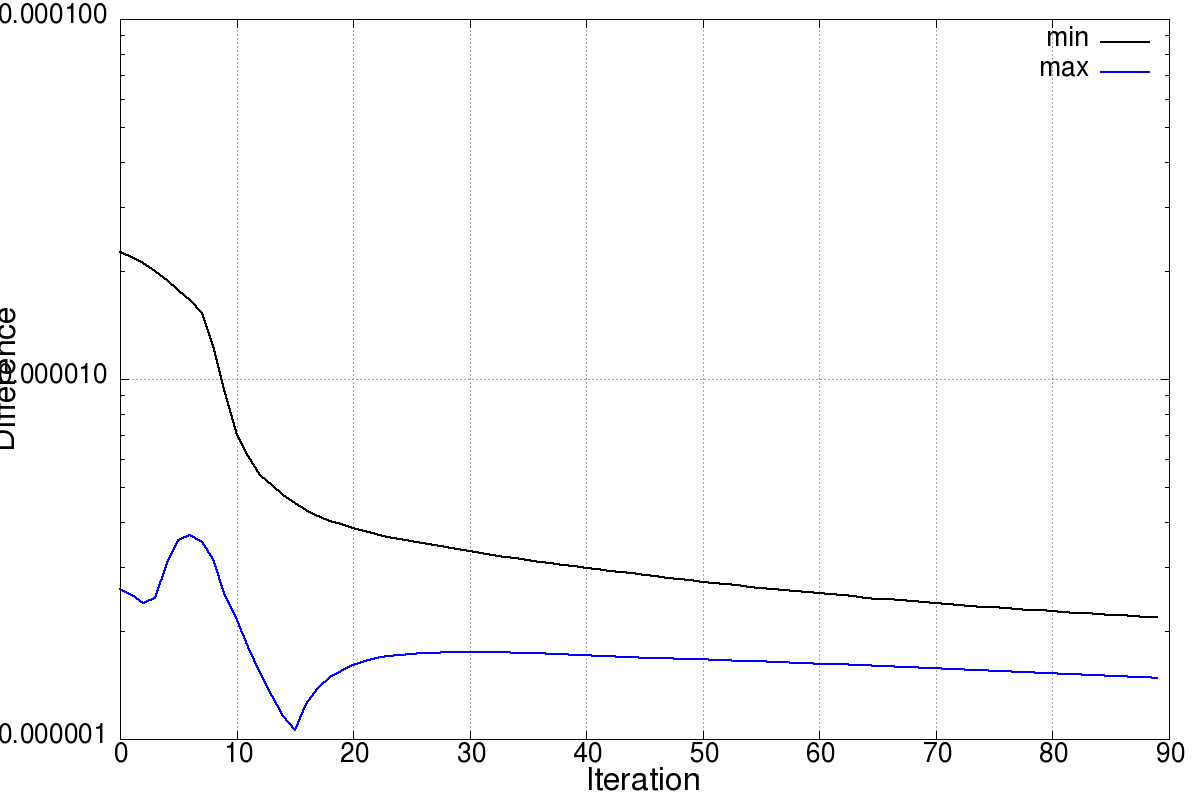}
    \caption{Plate With Hole: Start: $\alpha=1.0$, abs(Minimum) and Maximum}
    \label{fig:plate_with_hole_3}
\end{figure}

The case is shown in Figures \ref{fig:plate_with_hole_1}, 
\ref{fig:plate_with_hole_2} and \ref{fig:plate_with_hole_3} and
considers a plate with a hole. The plate dimensions are (all units in mks): 
$0 \le x \le 60$, $0 \le y \le 30$, $0 \le z \le 0.1$. A hole 
of diameter $d=10$ is placed in the middle ($x=30, y=15$).
Density, Young's modulus and Poisson rate were set to
$\rho=7,800, E=2 \cdot 10^{11}, \nu=0.3$ respectively.
672 linear, triangular, plain stress elements were used. 
The left boundary of the plate is assumed clamped ($\uvec=0$),
while a horizontal load of $q_x=10^5$ was prescribed at the 
right end.
The left part of the figures show the computed strength factor 
and displacements, while the right part displays the expected
values (the strength factor range is $0.1 \le \alpha \le 1$). 
The 14~measurement points, together with the differences
in displacements between measured and computed values are also
shown in the bottom right part.

\subsection{L-Shape}

\begin{figure}[!hbt]
    \centering
    \subfigure[Target displacements and sensor locations]{\includegraphics[width=0.35\linewidth]{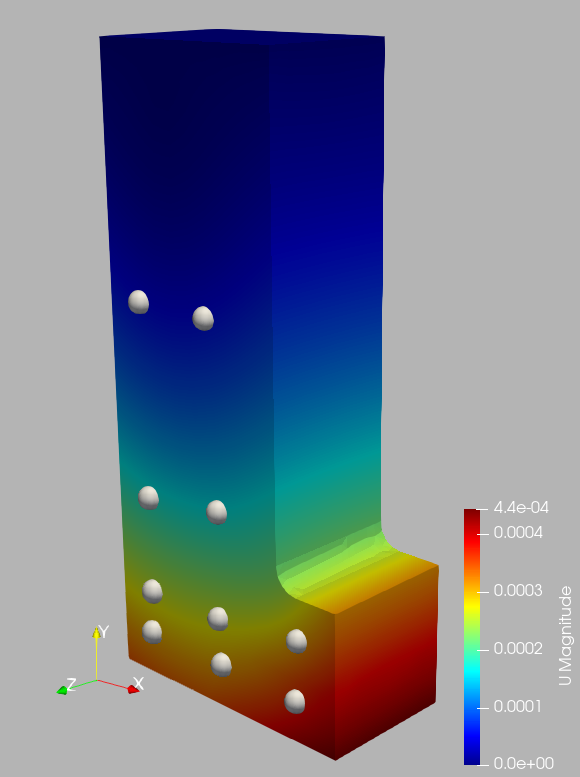}}
    \subfigure[Target strength factor]{\includegraphics[width=0.35\linewidth]{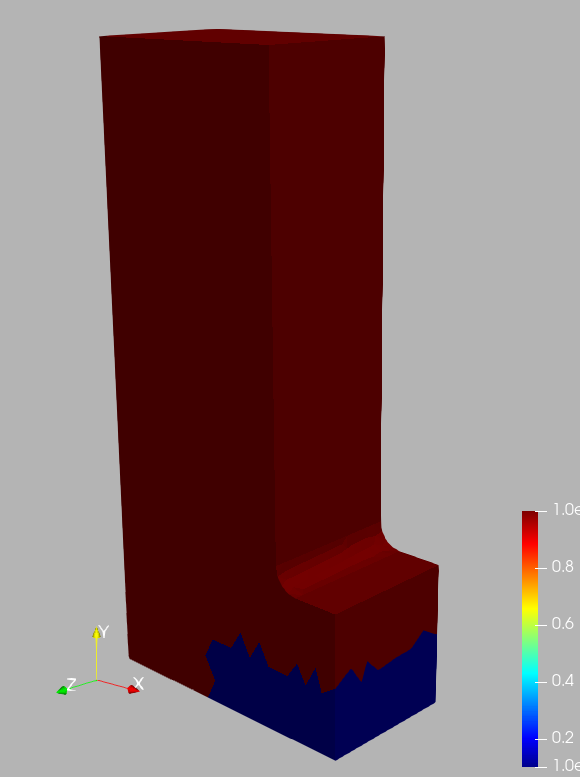}}
    \caption{L-Shape: Target conditions}
    \label{fig:lshape_target}
\end{figure}

\begin{figure}[!hbt]
    \centering
    \subfigure[Initial displacements and sensor locations]{\includegraphics[width=0.35\linewidth]{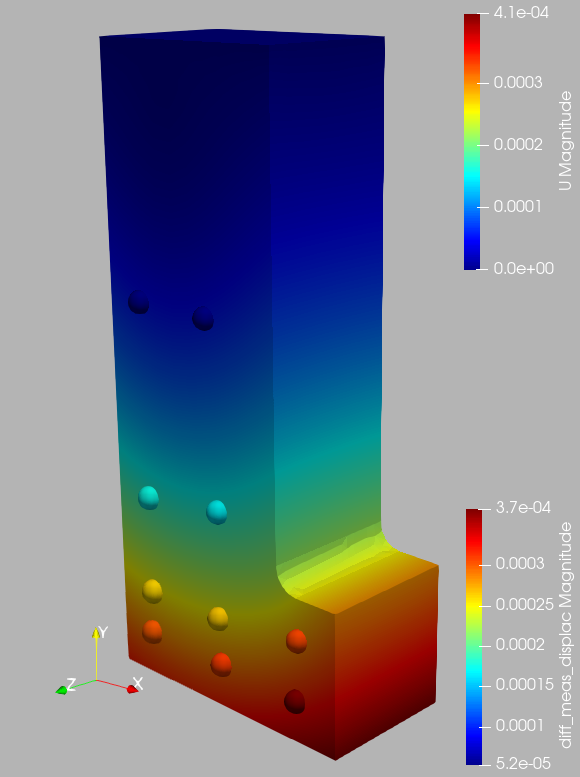}}
    \subfigure[Initial strength factor]{\includegraphics[width=0.35\linewidth]{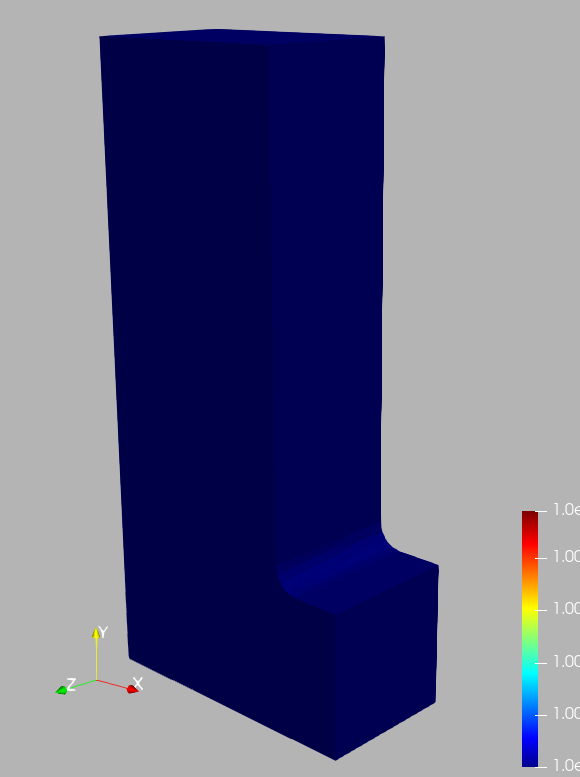}}
    \caption{L-Shape: Initial conditions}
    \label{fig:lshape_initial}
\end{figure}

\begin{figure}[!hbt]
    \centering
    \subfigure[Displacements obtained]{\includegraphics[width=0.35\linewidth]{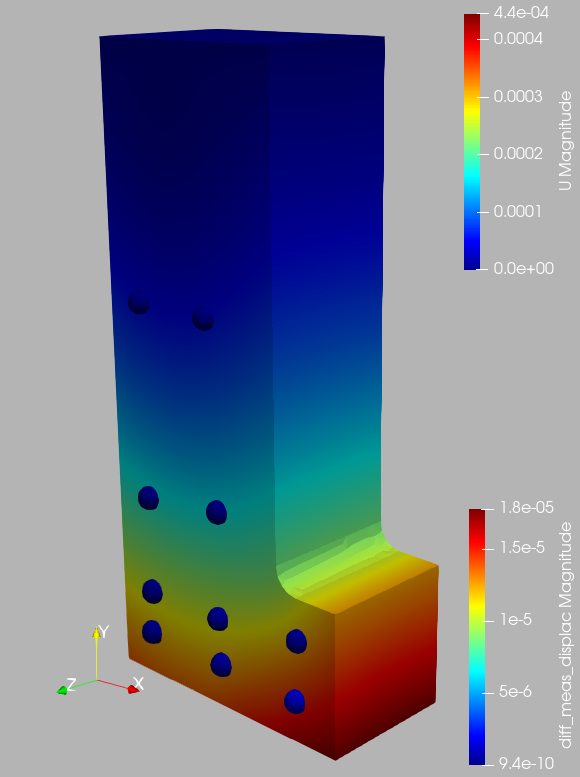}}
    \subfigure[Strength factor obtained]{\includegraphics[width=0.35\linewidth]{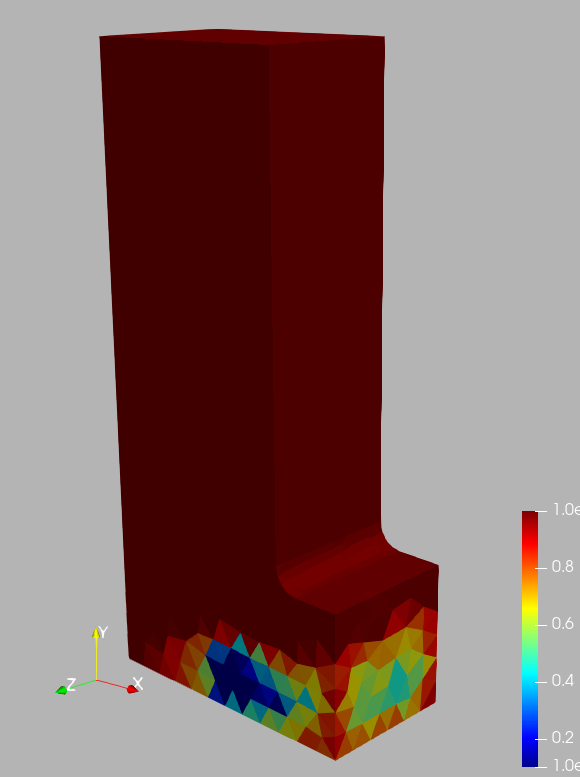}}
    \caption{L-Shape: Solution obtained}
    \label{fig:lshape_optim}
\end{figure}

The case, taken from \cite{di2022data} is shown in Figures 
\ref{fig:lshape_target}, \ref{fig:lshape_initial} and 
\ref{fig:lshape_optim} and considers 
an L-shaped block subjected to a vertical force.
The plate dimensions are (all units in mks):
$0 \le x \le 0.6$, $0 \le y \le 1.3$, $0 \le z \le 0.30$. The upper
part extends up to $x=0.45$, and the L-part extends to $y=0.3$.
A fillet with radius $r=0.05$ was added to avoid extreme stress
concentrations. Density, Young's modulus and Poisson rate were set to
$\rho=7,800, E=2 \cdot 10^{11}, \nu=0.3$ respectively.
14,622 linear, tetrahedral elements were used. 
The top boundary of the block is assumed clamped ($\uvec=0$),
while a vertical surface load of $f_y=-2 \cdot 10^7$ was prescribed 
at the top of the L-shaped region (only the straight section, i.e. not
the fillet).

The 10 visible measurement points (the other
10 are at the same $x,y$ positions but on the other $z$-face), together 
with the target displacements and strength factors are shown in 
figure \ref{fig:lshape_target} (the strength factor range again is 
$0.1 \le \alpha \le 1$).
. This case was particularly challenging 
because the weakened region does not have a considerable influence on the 
displacements. Therefore, many possible strength factor 
distributions can yield similar displacements. The smoothing of the gradient 
was a key tool for the optimizer to arrive at the proper solution.

\begin{figure}
    \centering
    \includegraphics[width=0.8\textwidth]{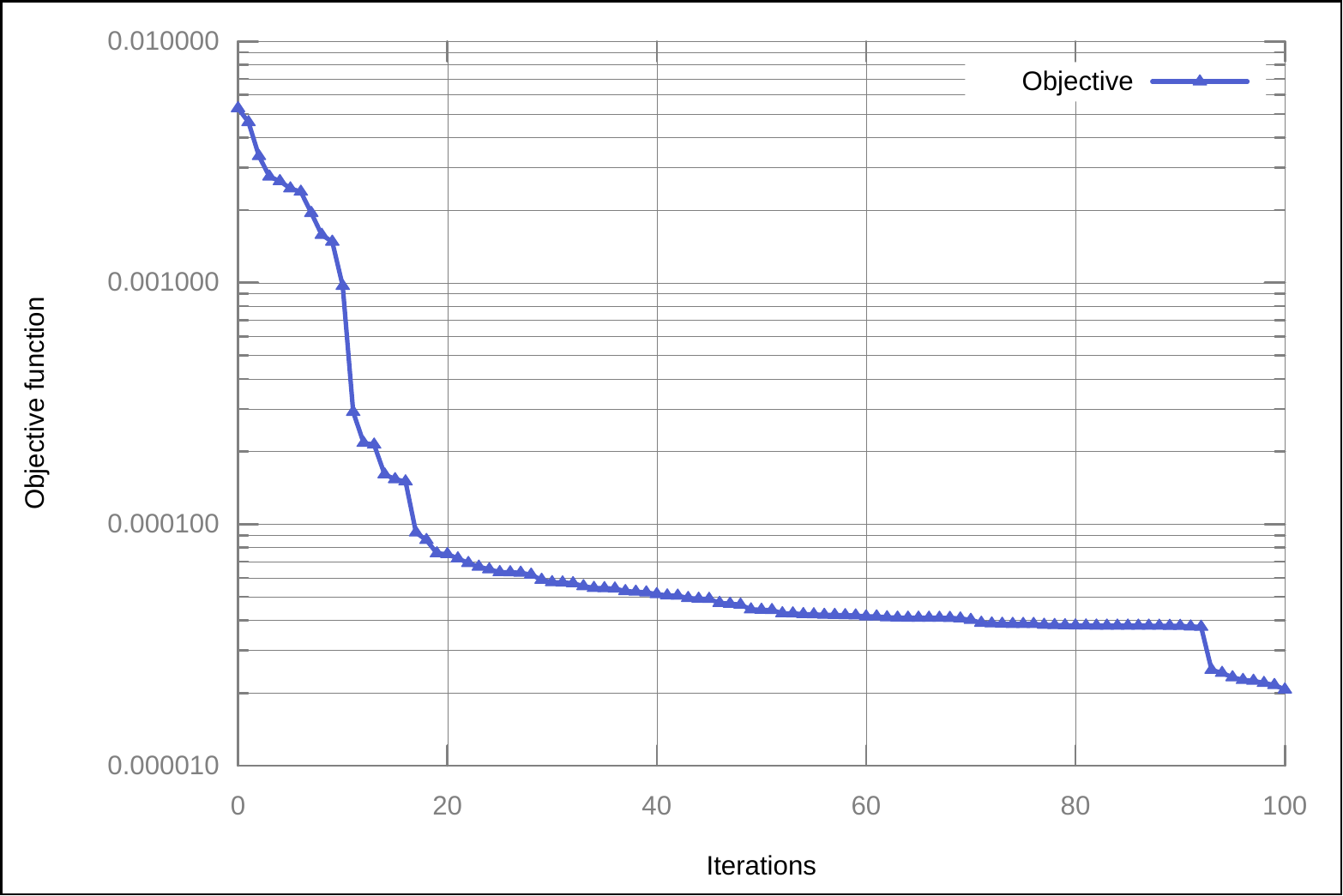}
    \caption{L-Shape: Objective function history}
    \label{fig:lshape_obj}
\end{figure}

Figure \ref{fig:lshape_initial} shows the initial conditions for the 
optimization loop. The results obtained after 100 steepest descent 
iterations are displayed in figure \ref{fig:lshape_optim}. Note that 5 passes 
of gradient smoothing were employed. The evolution of the objective function 
is shown in figure \ref{fig:lshape_obj}.

\section{Conclusions and outlook} \label{sec:conclusions}

An adjoint-based procedure to determine weaknesses, or, more generally
the material properties of structures has been presented. Given a
series of force and deformation/strain measurements, the material
properties are obtained by minimizing the weighted differences 
between the measured and computed values. It was found that in order
to obtain reliable, convergent results the gradient of the cost 
function has to be smoothed.
\par \noi
Several examples are included that show the viability, accuracy 
and efficiency of the proposed methodology using both displacement
and strain measurements.
\par \noi
We consider this a first step that demonstrates the viability of the
adjoint-based methodology for system identification and its use for
digital twins \cite{mainini2015surrogate, chinesta2020virtual}. 
Many questions remain open, of which we just mention two obvious ones:
\par \noi
\begin{itemize}
\item[-] What sensor resolution is required to obtain reliable results ?
\item[-] Will these techniques work under uncertain measurements ? 
\cite{HAntil_SDolgov_AOnwunta_2022b,HAntil_SDolgov_AOnwunta_2023a}.
\end{itemize}
\par \noi
Furthermore, the steepest descent procedures may be improved
by going to a quasi or full Newton solver. But: will they be
faster ?
\par \noi
The answers to these questions are currently under investigation.

\section{Acknowledgements}

This work is partially supported by NSF grant DMS-2110263 and the 
AirForce Office of Scientific Research under Award NO: FA9550-22-1-0248.

\bibliographystyle{plain}
\bibliography{references.bib}  

\end{document}